\documentclass{cocv}
\begin{document}
%%-----------------------------
%%      the top matter
%%-----------------------------
%\input{tcilatex}

\title{Boundary null-controllability of two coupled parabolic equations : simultaneous condensation of eigenvalues and eigenfunctions}

\author{El hadji Samb }\address{Aix Marseille Universit\'e, CNRS, Centrale Marseille, I2M, UMR 7373, 13453 Marseille, France.}

\thanks{The author thanks A. Benabdallah, F. Boyer et M. Morancey for useful conversations. I am very grateful to them.}

%\date{...}
\date\today

\begin{abstract} 
Let the matrix operator $L=D\partial_{xx}+q(x)A_0 $, with  $D=diag(1,\nu)$, $\nu\neq 1$, $q\in L^{\infty}(0,\pi)$, and $A_0$ is a Jordan block of order $1$. We analyze the boundary null controllability  for the system $y_{t}-Ly=0$. When $\sqrt{\nu} \notin \mathbb{Q}_{+}^*$ and  $q$ is constant, $q=1$ for instance, there exists a family of root vectors of $(L^*,\mathcal{D}(L^*))$ forming a Riesz basis of $L^{2}(0,\pi;\mathbb{R}^2 )$. Moreover in  \cite{JFA14} the authors show the existence of a minimal time of control depending on condensation of eigenvalues of $(L^*,\mathcal{D}(L^*))$, that is to say the existence of $T_0(\nu)$ such that the system is null controllable at time $T > T_0(\nu)$ and not null controllable at time  $T < T_0(\nu)$. In the same paper, the authors prove that for all $\tau \in [0, +\infty]$, there exists $\nu \in ]0, +\infty[$ such that $T_0(\nu)=\tau$. When $q$ depends on $x$, the property of Riesz basis is no more guaranteed. This leads to a new phenomena: simultaneous condensation of eigenvalues and eigenfunctions. This condensation affects the time of null controllability.
\end{abstract}

\maketitle

%%%%%%%%%%%%%%%%%%%%%%%%%%%%%%%%% introduction%%%%%%%%%%%%%%%%%
\section{Introduction and main results }
%%%%%%%%%%%%%%%%%%%%%%%%%%%%%%%%%%%%%%%%%%%%%%%%%%

This paper deals with the controllability of two coupled one-dimensional
parabolic equations, with different diffusion coefficients, where the control is exerted at one boundary point.
Let us fix $T > 0$ and consider the following control problem : 

\begin{equation} 
\left\{ 
	\begin{aligned}
& y_{t}+Ly=0 & \hbox{ in }&Q_{T}:= (0,\pi )\times
(0,T), \\ 
& y(0,\cdot )=Bu,\quad y(\pi ,\cdot )=0 & \hbox{ on }&(0,T),\\
& y(\cdot ,0)=y_{0} &\hbox{ in }&(0,\pi ),%
\end{aligned}%
\right.  
\label{boun.sg2x2}
\end{equation}%
where
\begin{equation} 
%\begin{aligned}
L=-D\partial_{xx}+q(x)A_0, \:
D = \hbox{diag} (1, \nu),\:\hbox{with $\nu > 0$}, \:  A_0=\left( 
\begin{array}{cc}
0 & 1 \\ 
0 & 0%
\end{array}%
\right) \hbox{ and } \:
B=\left( 
\begin{array}{c}
0 \\ 
1%
\end{array}%
\right),  
%\end{aligned}
\label{matDAB}
\end{equation} 
$q\in L^{\infty}(0,\pi)$ is a given function, $y_0$ is the initial datum  and $u \in L^2(0,T )$ is the control function. Equivalently, the previous system (\ref{boun.sg2x2}) can be written as

\begin{equation} 
\left\{ 
	\begin{aligned}
&	\partial_{t}y_{1}=\partial_{xx}y_1 - qy_2 & \hbox{ in }&Q_{T}:= (0,\pi )\times
(0,T), \\ 
&\partial_{t}y_{2}=\nu \partial_{xx}y_2 & \hbox{ in }&Q_{T}, \\ 
&	y_1(0,\cdot )=0,\quad y_2(0,\cdot )=u,\quad y(\pi ,\cdot )=0 & \hbox{ on }&(0,T),\\
&y(\cdot ,0)=y_{0} &\hbox{ in }&(0,\pi ).%
\end{aligned}%
\right. 
\label{boun.sg2x2+}
\end{equation}
It is known (see \cite{FGT} Prop. 2.2) that for any given initial data $y_{0}\in H^{-1}(0,\pi;\mathbb{R}^2 )$  and  $u\in L^{2}(0,T)$, system \eqref{boun.sg2x2} possesses a unique solution defined by transposition which satisfies 
$$y\in L^{2}\left(Q_T;\mathbb{R}^2\right)\cap C^{0}\left([0,T];H^{-1}(0,\pi;\mathbb{R}^2)\right)$$
and depends continuously on the data $u$ and $y_0$, i.e., there exists a constant $C = C(T) > 0$ such that

\begin{equation*} 
\Vert y\Vert_{ L^{2}\left(Q_T;\mathbb{R}^2\right)}+\Vert y\Vert_{C^{0}\left([0,T];H^{-1}(0,\pi;\mathbb{R}^2)\right)}\leq C\left( \Vert y_0\Vert_{ H^{-1}(0,\pi;\mathbb{R}^2)}+\Vert u\Vert_{ L^{2}\left(0,T\right)} \right).
%\label{solbounb en}
\end{equation*}
Let us introduce the notion of null and approximate controllability for this kind of systems.

%\begin{dfntn}
\begin{enumerate}

\item System (\ref{boun.sg2x2}) is approximately controllable in $H^{-1}(0,\pi;\mathbb{R}^2 )$ at time $T$ if for every  $y_0,y_d \in H^{-1}(0,\pi;\mathbb{R}^2 )$ and for every $\varepsilon>0$, there exists a control $u\in L^{2}(0,T)$ such that the solution to  (\ref{boun.sg2x2}) associated to $y_0$ and $u$ satisfies 
$$
\Vert y(\cdot,T)-y_d \Vert _{H^{-1}(0,\pi;\mathbb{R}^2 )} \le \varepsilon.
$$
\item 
System (\ref{boun.sg2x2}) is null controllable at time $T$ if for every initial
condition $y_0\in H^{-1}(0,\pi;\mathbb{R}^2 )$, there exists a control $u\in L^{2}(0,T)$ such that the
solution $y$ to (\ref{boun.sg2x2}) satisfies 
$$
y(\cdot, T) =0 \quad \hbox{in}\; H^{-1}(0,\pi;\mathbb{R}^2).
$$
\end{enumerate}
%\end{dfntn}

The null controllability of parabolic partial differential equations has been widely studied since the pioneering work of \cite{FATRUS}. From the works of  \cite{FI96} and  \cite{LR95}, it was commonly admitted that, in the context of parabolic partial differential equations, there is no restriction on the final time $T$. But recently the study of particular examples highlighted the existence of a positive \textit{minimal time} $T_0$ for null controllability. Actually, such an example was already provided in the 70s in \cite{DOL73}. The more recent results concerning such strictly positive a minimal time have been obtained in contexts of control of coupled parabolic equations \cite{JFA14}, \cite{JMAA16}, \cite{DZesaimcocv17} and more generally see \cite{JMPA17}.

The main goal of this article is to address a new phenomenon arising in the null controllability issue for system (\ref{boun.sg2x2}). Let us first recall some known results about the controllability properties of scalar parabolic systems. The null controllability problem for scalar
parabolic systems has been first considered in the one-dimensional case. In \cite{FR1} and \cite{FATRUS}, using the \textit{moment method}, H. O. Fattorini and D. L. Russell gave a positive answer for null controllability problem of considered parabolic system.  Also, the authors proved a general result on existence of a bi-orthogonal family to $\{ e^{-\Lambda_k t}\}_{k \ge 1}$ in $ L^2(0, T)$ which fulfils appropriate bounds if the sequence $\Lambda_k \subset \mathbb{R}_+$ satisfies the following \textit{gap property}:
\begin{equation}
\sum_{k \ge 1}\frac{1}{\Lambda_k} < \infty \quad and \quad \vert  \Lambda_k-\Lambda_l\vert \ge \rho\vert k-l\vert , \quad \forall k,l \ge 1,
\label{gapprop}
\end{equation}%
for a constant $\rho >0$. In 1973, S. Dolecki addressed the pointwise controllability at time $T$ of the one-dimensional heat equation (see \cite{DOL73}). S. Dolecki exhibited a minimal time  which depends on the point $x_0$. To our knowledge, this was the first result on null-controllability of parabolic problems where a minimal time of control appears. System (\ref{boun.sg2x2}) is a particular class of more general $n\times n$ parabolic control systems
of the form:

\begin{equation} 
\left\{ 
	\begin{array}{ll}
\partial_ty-(D\partial_{xx}+A(x))y=0 & \hbox{in }Q_{T}:= (0,\pi )\times
(0,T), \\ 
	\noalign{\smallskip}
y(0,\cdot )=Bu,\quad y(\pi ,\cdot )=0 & \hbox{ on }(0,T),\\
y(\cdot ,0)=y_{0} &\hbox{ in }(0,\pi ).%
\end{array}%
\right.  \label{boun.sg2x2gene}
\end{equation}%
Here, $D=diag(d_1,..,d_n)$, with $d_i>0$ for $i:1\le i\le n$, $A=(a_{ij})_{1\le i\le n}\in L^{\infty}((0,\pi);M_n(\mathbb{R}))$ and $B \in \mathcal{L}(\mathbb{R}^m,\mathbb{R}^n) $. In system (\ref{boun.sg2x2gene}), $u\in L^2(0,T;\mathbb{R}^m) $ is the control and we want to control the complete system ($n$ equations) by means of $m$ controls exerted on the boundary condition at point $x = 0$. Observe that the most interesting (and difficult) case is the case $m < n$.

The first results of null controllability for system (\ref{boun.sg2x2gene}) was obtained in \cite{FGT} in the case $n = 2$, $m = 1$ and $D = Id$, $A\in M_n(\mathbb{R})$. This result was generalized by \cite{JMPA11} to the case $n \ge 2$, $m \ge 1$.
In these two papers, the authors used the method of moments of Fattorini-Russell to give a necessary and sufficient condition of null controllability at any time $T > 0$ for system (\ref{boun.sg2x2gene}). In both cases, the  sequence of eigenvalues $\Lambda=\{  \Lambda_k\}_{k \ge 1}\subset \mathbb{R}_+$ of the matrix operator $\mathcal{A}=Id\partial_{xx}+A$ with Dirichlet boundary conditions continue to satisfy the \textit{gap} condition in (\ref{gapprop}). As in the scalar case, this gap property (together with
appropriate properties for the coupling and control matrices $A$ and $B$) provides the null controllability result for system (\ref{boun.sg2x2gene}) at any positive time.

In \cite{JFA14}, the authors are interested in the extension of the previous null controllability  results for system (\ref{boun.sg2x2gene}) to the case where $D \neq Id$, $A\in M_n(\mathbb{R})$, $n > 1$ and $m < n$. The main difference with the case $D = Id$ lies in the behavior of the sequence of eigenvalues of the matrix operator $\mathcal{A} := D\partial_{xx} + A$. The operator $-\mathcal{A}$ admits a sequence of eigenvalues $\Lambda=\{ \Lambda_k \}_{k \ge 1}$  which does not satisfy the gap condition appearing in (\ref{gapprop}) but the operator $-\mathcal{A}$ is  diagonalizable, i.e., its eigenfunctions form a Riesz basis. As a consequence, the authors show the existence of a minimal time of control $T_0 \in [0,+\infty]$,  depending to the so-called \textit{condensation index}, $c(\Lambda)$, of the sequence $\Lambda$ of eigenvalues of the operator $-\mathcal{A}$. To our knowledge, this condensation index has been introduced for the first time by V. Bernstein (see \cite{BER}) for increasing real sequences and later extended by J. R. Shackell (see \cite{SHA}) to complex sequences. Roughly speaking, if we consider the complex sequence  $\Lambda=\{ \Lambda_k \}_{k \ge 1}\subset \mathbb{C}$, the condensation index of $\Lambda$,  is a measure of the way how $\Lambda_n$ approaches $\Lambda_m$ for $n \neq m$.

In \cite{JMAA16}, the authors consider the case where $D = Id$, $A(x)=q(x)A_0$, with $q\in L^{\infty}(0,\pi)$ , $n =2$ and $m =1$. In this case the  operator $L^* := -Id\partial_{xx} + q(x)A^*_0$ admits a sequence of eigenvalues $\Lambda=\{ \Lambda_k \}_{k \ge 1}$  which satisfy the gap condition appearing in (\ref{gapprop}) moreover using the eigenfunctions and the generalized eigenfunctions of the operator $L^*$, we can construct a Riesz basis for the space $L^{2}(0,\pi;\mathbb{R}^2)$ .  As a consequence, the authors show the existence of a minimal time of control $T_0(q)$. They also proved that for any $\tau \in [0,+\infty]$, there exists $q\in L^{\infty}(0,\pi)$ such that $T_0(q)=\tau$.

In this paper, we study the null controllability properties of system (\ref{boun.sg2x2gene}) to the case where $D \neq Id$ and $A(x)=q(x)A_0$, $q\in L^{\infty}(0,\pi)$,  $n =2$ and $m =1$. This situation may seem like a simple perturbation of previous cases ( in \cite{JFA14} and \cite{JMAA16}). It is not true, it contains a new phenomenon : simultaneous condensation of eigenvalues and eigenfunctions. This phenomenon was excluded from all previous cases because of the property that the family of eigenfunctions of the operator $(L^*,\mathcal{D}(L^*))$ form a Riesz basis for the  Hilbert space where the system is posed. This condensation of eigenfunctions can compensate the condensation of eigenvalues  and the  minimal time of control is affected. In \cite{JFA14} as in \cite{JMAA16} the authors proved the controllability by using the usual moments method but this method does not take into account such phenomena  of eigenfunctions condensation. Under appropriate assumptions on $D$ and $q\in L^{\infty}(0,\pi)$, the operator $L^*= -D\partial_{xx} +q(x)A^*_0$ admits a sequence of eigenvalues $\Lambda=\{ \Lambda_k \}_{k \ge 1}$  which does not satisfy the gap condition appearing in (\ref{gapprop}) however the sequence of associated eigenfunctions is \textit{complete} but it is not a Riesz basis (see Proposition \ref{notriesz}) for some $\sqrt{\nu} \notin \mathbb{Q}^*_+$ and $q\in L^{\infty}(0,\pi)$ . As a consequence, we will see that a minimal time of control $T_0 \in [0,+\infty]$,  depends simultaneously of condensation of eigenvalues and   associated eigenfunctions  of $(L^*,\mathcal{D}(L^*))$. To this end, we will use the \textit{block methods moment} developed by  A. Benabdallah, F. Boyer, M. Morancey ( in \cite{BBM18}).

The plan of the paper is the following one: In Section 1, we address some known results about the controllability of parabolic system and we give the main result of this work. In Section 2 we study the null-controllability of system (\ref{boun.sg2x2}) when $\sqrt{\nu}\notin \mathbb{Q}_+^*$, $\sqrt{\nu}>1 $ (resp. $\sqrt{\nu}<1 $). Section 3 is devoted to  null-controllability of system (\ref{boun.sg2x2}) when $\sqrt{\nu}\in \mathbb{Q}_+^*$. Finally, in the Appendix we give additional properties of our main result.

Let us present our first boundary control results, when  $\sqrt{\nu} \notin \mathbb{Q}_{+}^*$. Let us first introduce some notations. Let  $B$, $A_0$ and $D$ given by (\ref{matDAB}) and a function $q\in L^{\infty}(0,\pi)$. 
Let us consider the operator

\begin{equation}  \label{opL}
 L^*:=-D\partial_{xx}+q(x)A^{\ast}_0: \mathcal{D}(L^*) \subset L^{2}(0,\pi ;\mathbb{R}^{2}) \longrightarrow  L^{2}(0,\pi;\mathbb{R}^{2}) 
\end{equation}  
with domain $\mathcal{D}(L^*)=H^{2}(0,\pi ;\mathbb{R}^{2})\cap H_{0}^{1}(0,\pi ;\mathbb{R}^{2})$. Given $k\geq1$, let us consider
 
\begin{equation} \Phi _{1,k}^{\ast }:=\left( 
\begin{array}{l}
\varphi _{k} \\ 
\psi _{k}%
\end{array}%
\right) ,\Phi _{2,k}^{\ast }:=\left( 
\begin{array}{l}
0 \\ 
\varphi _{k}%
\end{array}%
\right),
 \label{vectpropnuir}
\end{equation}%
where $\psi_k$ is the unique solution of problem:

\begin{equation} 
\left\{ 
	\begin{array}{ll}
\nu \psi_k''+ k^2\psi_k=q\varphi_k \: & \hbox{in }(0,\pi ), \\ 
	\noalign{\smallskip}
\psi(0)=\psi(\pi )=0. 
\end{array}%
\right.
 \label{phikprob}
\end{equation}%
Assume that  $B^{\ast}D\frac{\partial\Phi_{1,k}^{\ast }}{\partial x} (0)\neq 0$, for all $k\ge 1$. Let us defined, for any $k \in \mathbb{N}^*$, $i_k$ (resp. $j_k$) as the nearest integer to $\sqrt{\nu}k$ (resp. $\frac{k}{\sqrt{\nu}}$ ), i.e,
 $ \left \vert \sqrt{\nu}k -i_k   \right \vert<\frac{1}{2}$ (resp. $\left \vert \frac{k}{\sqrt{\nu}} -j_k   \right \vert<\frac{1}{2}$ ). Let us denote
\begin{equation} \label{psik}
\psi_{i,k}:=\frac{\Phi_{i,k}^{\ast }}{B^{\ast}D\frac{\partial\Phi_{i,k}^{\ast }}{\partial x} (0)},\;\forall k\ge 1,\quad \forall i=1,2.
\end{equation}
In the same way, denote
\begin{equation} \label{psik+}
\psi_{1,i_k}:=\frac{\Phi_{1,i_k}^{\ast }}{B^{\ast}D\frac{\partial\Phi_{1,i_k}^{\ast }}{\partial x} (0)} \;  \hbox{ and }\; \psi_{2,j_k}:=\frac{\Phi_{2,j_k}^{\ast }}{B^{\ast}D\frac{\partial\Phi_{2,j_k}^{\ast }}{\partial x} (0)},\quad \forall k\ge 1.
\end{equation}
One has
%%%%%%%%%%%%%%%%%%%%%%Thm1%%%%%%%%%%%%%%%%%
\begin{thrm} \label{Thm}
Suppose $\sqrt{\nu} \notin \mathbb{Q}_{+}^*$. 

\begin{enumerate}

\item The spectrum of $\left( L^{\ast },\mathcal{D}(L^*)\right)$ is given by $\sigma(L^*)=\{k^2\}_{k\geq 1} \cup \{\nu k^2\}_{k\geq 1}$  and the corresponding family of eigenfunctions, $\left\lbrace \Phi_{1,k}^{\ast },\Phi_{2,k}^{\ast } \right\rbrace_{k\geq 1}$, is complete in $L^{2}(0,\pi;\mathbb{R}^2 )$. %but is not a Riesz basis for $L^{2}(0,\pi;\mathbb{R}^2 )$.

\item System (\ref{boun.sg2x2}) is approximately controllable at time $T > 0$, if and only if
\begin{equation} \label{cndca}
B^{\ast}D\frac{\partial\Phi_{1,k}^{\ast }}{\partial x} (0)\neq 0,\;\forall k\ge 1.
\end{equation}
\item

Introduce
\begin{equation}  \label{T0tilde}
\begin{aligned}
 \widetilde{T}_0  :=  \quad \max \left\lbrace \displaystyle{\limsup_{k \to +\infty } } \frac{\log \left \Vert \psi_{1,k}\right \Vert_{H^1_0}}{k^2}, \displaystyle{\limsup_{k \to +\infty } }  \frac{\log \frac{\left \Vert \psi_{1,i_k}-\psi_{2,k} \right \Vert_{H^1_0}}{\vert  i_k^2-\nu k^2\vert }}{ \nu k^2}\right\rbrace=\max \left\lbrace T_1,  \widetilde{T}_2 \right\rbrace \, \in [0, +\infty], 
 \end{aligned}
 \end{equation}
 and 
 \begin{equation}  \label{T0hat}
 \begin{aligned}
 \widehat{T}_0 & :=  \quad \max \left\lbrace \displaystyle{\limsup_{k \to +\infty } } \frac{\log \left \Vert \psi_{1,k}\right \Vert_{H^1_0}}{k^2},\displaystyle{\limsup_{k \to +\infty } } \frac{\log \frac{\left \Vert \psi_{1,k}-\psi_{2,j_k} \right \Vert_{H^1_0}}{\vert  k^2-\nu j_k^2\vert }}{ k^2}\right\rbrace =  \max \left\lbrace  T_1,  \widehat{T}_2 \right\rbrace \, \in [0, +\infty].
  \end{aligned}
\end{equation}

\begin{enumerate}
\item 
 If $\sqrt{\nu}>1$ then : system~(\ref{boun.sg2x2}) is null-controllable in $H^{-1}(0,\pi;\mathbb{R}^2)$  for  $T> \widetilde{T}_0$. On the other hand, if $\widetilde{T}_0>0$, system~(\ref{boun.sg2x2}) is not null-controllable in $H^{-1}(0,\pi;\mathbb{R}^2)$  for  $T< \widetilde{T}_0$.

\item
 If $\sqrt{\nu}<1$ then : system~(\ref{boun.sg2x2}) is null-controllable in $H^{-1}(0,\pi;\mathbb{R}^2)$  for  $T> \widehat{T}_0$. On the other hand, if $\widehat{T}_0>0$, system~(\ref{boun.sg2x2}) is not null-controllable in $H^{-1}(0,\pi;\mathbb{R}^2)$  for  $T< \widehat{T}_0$.
\end{enumerate}
\end{enumerate}
\end{thrm}
%%%%%%%%%%%%%%%%%%%%%%%%%%%%%%%%%%%%%%%%%%%Thm%%%%%%%%%%%%%%%%%%%%%%%%%

\smallskip

\begin{rmrk}
\label{remark1thm} 

\begin{enumerate}
\item In Appendix \ref{AppA}, we will write in a simple form the expressions of Theorem \ref{Thm} and  show that, $\widetilde{T}_0$ and $\widehat{T}_0$ strongly depends on diffusion coefficient $\sqrt{\nu} \notin \mathbb{Q}_{+}^*$ and the coupling function $q\in L^{\infty}(0,\pi)$. Moreover we will remark that $\widehat{T}_0=\widehat{T}_2\ge T_1 $ and $\widehat{T}_2\ge\widetilde{T}_2 $ (see (\ref{T0tilde}) and (\ref{T0hat})).

\item Condition (\ref{cndca}) characterizes the approximate controllability property of system (\ref{boun.sg2x2}). Thus, (\ref{cndca}) is a necessary
condition for the null controllability at time $T > 0$ of this system.

    \item 
The approximate controllability result stated in Theorem \ref{Thm} does not depend on the final time $T$: approximate controllability of system (\ref{boun.sg2x2}) at a time $T_0 > 0$ is equivalent to the
approximate controllability of system  at any time $T > 0$.

\item We will prove, in  Proposition \ref{prop exmple}, that $\widetilde{T }_0$ (resp. $\widehat{T }_0$)  could take any value on the interval $[0,\infty]$. Thus, when $\widetilde{T }_0,\widehat{T }_0\in(0,\infty]$, from Theorem \ref{Thm}
we deduce that system (\ref{boun.sg2x2}) could have the approximate controllability property at a positive time $T$ without being null controllable at this time $T$.

%\item
%In particular, Theorem \ref{Thm} implies that, even
%in a parabolic setting, a positive time of control may appear and that, unlike the scalar case, boundary approximate and null controllability are not equivalent properties in the non-scalar case. For such functions, system (\ref{boun.sg2x2}) is approximately controllable at all positive time $T$ but it is not null controllable at time $T$ if $T \in (0, \widetilde{T}_0)$ (or if $T \in (0, \widehat{T}_0)$).
\end{enumerate}
\end{rmrk} 

Let us now present our second boundary control results when  $\sqrt{\nu}=\frac{i_0}{j_0} \in \mathbb{Q}_{+}^*$ , where $i_0 $ and $j_0 $ are co-prime ( $i_0 \land j_0=1$ ). 

%%%%%%%%%%%%%%%%%%%%%%%%%%%%%%%%%%%%%%%%%%%%%%%%%%%%%%Thm+%%%%%%%%%%%%%%%%%%%%%%%%%%%
\begin{thrm} \label{Thm+}
Let us consider $\sqrt{\nu}\in \mathbb{Q}^*_{+}$. Let  $B$, $A_0$ and $D$ given by (\ref{matDAB}).  To the function $q\in L^{\infty}(0,\pi)$ and $\nu$, we associate the sequence $\{I_{k}\}_{k\ge 1}$ , defined by 
\begin{equation}
 I_{k}(\nu,q)=\int_{0}^{\pi}q(s)\varphi_{k}(s)\sqrt{\frac{\pi}{2}}\sin\left(\frac{k}{\sqrt{\nu}}(\pi-s)\right)ds.
 \label{Ik-}
\end{equation}
Then, one has:
\begin{enumerate}
\item System (\ref{boun.sg2x2}) is approximately controllable at time $T > 0$, if and only if
\begin{equation} \label{cndca+}
 I_{k}(\nu,q)\neq 0,\quad\forall k\ge 1.
\end{equation}
\item Assume (\ref{cndca+}) holds and define 
\begin{equation}  \label{T0+}
  T_0(\nu,q):=\limsup_{k \to +\infty }  \frac{-\log  \left \vert I_{k}(\nu,q)  \right \vert }{ k^2} \in [0, +\infty].
\end{equation}
Then if $T>T_0(\nu,q)$ system~(\ref{boun.sg2x2}) is null-controllable at time $T$. On the other hand, if $T_0(\nu,q)>0$, system~(\ref{boun.sg2x2}) is not null-controllable at time $T$ in $H^{-1}(0,\pi;\mathbb{R}^2)$ for $T<T_0(\nu,q)$.
\end{enumerate}
\end{thrm}
%%%%%%%%%%%%%%%%%%%%%%%%%%%%%%%%%%%%%%%%%%%%Thm+%%%%%%%%%%%%%%%%%%%%%%%%%%%%%%%
\smallskip

\begin{rmrk}
Under condition (\ref{cndca+}), the minimal time $T_0(\nu,q)$ is well-defined. Moreover, the sequence $\{ I_{k}(\nu,q)\}_{k\ge1}$ is bounded and $T_0(\nu,q)\in [0, +\infty]$.
\end{rmrk}

Let $q \in L^{\infty}(0,\pi)$ and $A_0$ given by (\ref{matDAB}). We introduce the backward adjoint problem associated with system~(\ref{boun.sg2x2}):
	\begin{equation}
	\left\{ 
	\begin{array}{ll}
-\theta _{t}-(D\partial _{xx}-qA_0^{\ast })\theta=0  & \hbox{in } Q_{T} , \\ 
	\noalign{\smallskip}
\theta(0,\cdot )=\theta(\pi ,\cdot )=0 & \hbox{ on }(0,T),\\
\theta(\cdot ,T)=\theta_{0} &\hbox{ in }(0,\pi ),%
\end{array}%
\right.  \label{Padjoint}
	\end{equation}
where $\theta_0 \in L^{2}(0,\pi ;\mathbb{R}^{2})$ is a given initial datum. Let us first see that this problem is well posed.
One has:

\begin{prpstn}
\label{prop1}
Assume that  $\theta_0 \in L^{2}(0,\pi ;\mathbb{R}^{2})$ is given. Then, system (\ref{Padjoint}) admits a unique solution $\theta \in L^2(0,T; H^1_0(0,\pi; \mathbb{R}^2))\cap C^0([0,T];L^2(0,\pi; \mathbb{R}^2))$ and in addition satisfies 
$$
\left \Vert \theta \right \Vert_{ L^2(0,T; H^1_0(0,\pi; \mathbb{R}^2))}+ \left \Vert \theta \right \Vert_{ C^0([0,T];L^2(0,\pi; \mathbb{R}^2))} \le C(T) \left \Vert \theta_0 \right \Vert_{ L^2(0,\pi; \mathbb{R}^2)},
$$
for a positive constant $C(T)>0$ independent of $\theta_0$. Furthermore, if  $\theta_0 \in H^{1}_0(0,\pi ;\mathbb{R}^{2})$, then the solution $\theta$ of the adjoint problem (\ref{Padjoint}) satisfies
$$
\theta \in L^2(0,T; H^2(0,\pi; \mathbb{R}^2) \cap H^1_0(0,\pi; \mathbb{R}^2))\cap C^0([0,T];H^1_0(0,\pi; \mathbb{R}^2))
$$
and, for a positive constant $C(T)>0$
$$
\left \Vert \theta \right \Vert_{ L^2(0,T; H^2(0,\pi; \mathbb{R}^2) \cap H^1_0(0,\pi; \mathbb{R}^2))}+ \left \Vert \theta \right \Vert_{ C^0([0,T];H^1_0(0,\pi; \mathbb{R}^2))} \le C(T) \left \Vert \theta_0 \right \Vert_{ H^1_0(0,\pi; \mathbb{R}^2)}
$$
\end{prpstn}

\bigskip

The next proposition provides a relation between systems (\ref{boun.sg2x2}) and (\ref{Padjoint}):

\begin{prpstn}
\label{prop2}
Let us consider  $ A_0 $ and $B$ given by (\ref{matDAB}) and $ q\in L^{\infty}(0,\pi)$. Then, for any $ y_0 \in H^{-1}(0,\pi; \mathbb{R}^2)$, $ u \in L^{2}(0,T)$ and $\theta_0 \in H^{1}_0(0,\pi; \mathbb{R}^2)$, one has 

\begin{equation*}
%\label{IPP en}
\int_{0}^{T}  u(t)B^{\ast }D\theta_x(0,t) dt = \langle y(\cdot , T) , \theta _0 \rangle_{H^{-1},H_0^1} - \langle y_0 , \theta (\cdot , 0) \rangle_{H^{-1},H_0^1} ,
\end{equation*}

where  $y\in L^{2}\left(Q_T;\mathbb{R}^2\right)\cap C^{0}\left([0,T];H^{-1}(0,\pi;\mathbb{R}^2)\right)$ and $$\theta \in L^2(0,T; H^2(0,\pi; \mathbb{R}^2) \cap H^1_0(0,\pi; \mathbb{R}^2))\cap C^0([0,T];H^1_0(0,\pi; \mathbb{R}^2))$$ are, resp., the solution to (\ref{boun.sg2x2}) and (\ref{Padjoint}) associated with $(u, y_0)$ and $\theta_0$.

\end{prpstn}

For a proof of the previous results see for instance \cite{TUW} or \cite{FGT}.

\bigskip

The controllability of system (\ref{boun.sg2x2}) can be characterized in terms of appropriate properties of
the solutions to the adjoint problem (\ref{Padjoint}). More precisely, we have

\begin{prpstn}
\label{dual_ncob en}
The following properties are equivalent:
\begin{enumerate}
\item There exists a positive constant $C>0$ such that, for any $ y_0 \in H^{-1}(0,\pi; \mathbb{R}^2)$, there exists a control $ u \in L^{2}(0,T)$ such that
\begin{equation*}
%\label{ubound}
\Vert u\Vert_{L^{2}(0,T) }\le C \Vert y_0\Vert_{H^{-1}(0,\pi; \mathbb{R}^2)}
\end{equation*}

and the associated state satisfies 

\begin{equation*}
%\label{yT=0 en}
y(\cdot, T) =0 \quad \hbox{in} \;H^{-1}(0,\pi;\mathbb{R}^2 ).
\end{equation*}
\item
There exists a positive constant $C$ such that the observability inequality
\begin{equation}
\label{InOb en}
\left \Vert \theta (\cdot, 0) \right \Vert^2_{H^1_0(0,\pi; \mathbb{R}^2)} \le C \int_{0}^{T}  \left \vert B^{\ast }D\theta_x(0,t) \right \vert^2 dt
\end{equation}
holds for every  $\theta_0 \in H^{1}_0(0,\pi ;\mathbb{R}^{2})$. In (\ref{InOb en}), $\theta$ is the adjoint state associated with $\theta_0$.
\end{enumerate}
\end{prpstn} 
This result is well known, the "Hilbert Uniqueness Method", see \cite{Lions68}.

%%%%%%%%%%%%%%%%%%%%%%%%%%%%%%%%%%%%%%%%%%%%%%%%%%%%%
\section{ Approximate and null-controllability result with an irrational diffusion coefficient  }
%%%%%%%%%%%%%%%%%%%%%%%%%%%%%%%%%%%%%%%%%%%%%%%%%%%%%%%%%%

%%%%%%%%%%%%%%%%%%%%%%%%%%%%%%%%%%%%%%%
\subsection{ Some preliminary results}
%%%%%%%%%%%%%%%%%%%%%%%%%%%%%%%%%%%%%%%%%%%%%%%%%
Let us consider the vectorial operator
\begin{equation*} % \label{opL}
 L:=-D\partial_{xx}+q(x)A_0: \mathcal{D}(L) \subset L^{2}(0,\pi ;\mathbb{R}
^{2})\longrightarrow L^{2}(0,\pi ;\mathbb{R}^{2}) 
\end{equation*}  
with domain $\mathcal{D}(L)=H^{2}(0,\pi ;\mathbb{R}^{2})\cap H_{0}^{1}(0,\pi ;\mathbb{R}^{2})$ and also its adjoint $L^{\ast }$ given by (\ref{opL}). % We will always denote by $\langle \cdot,\cdot \rangle_{X',X}$ the duality pairing between the Hilert space $X$ and its dual $X'$.

\begin{prpstn}\label{fam_nu_ir}
Let $A_0$ be given by ~(\ref{matDAB}) and consider the operator $L$ given by ~(\ref{opL}) and also denote its adjoint $L^{\ast }$. Assume that $\sqrt{\nu} \notin \mathbb{Q}$, then, 
\begin{enumerate}
\item The spectrum of $L^{\ast }$ is given by $\sigma(L^{\ast })=\{k^2\}_{k\geq 1} \cup \{\nu k^2\}_{k\geq 1}.$
\item Given $k\geq1$, if
 
\begin{equation*} \Phi _{1,k}^{\ast }:=\left( 
\begin{array}{l}
\varphi _{k} \\ 
\psi _{k}%
\end{array}%
\right) ,\Phi _{2,k}^{\ast }:=\left( 
\begin{array}{l}
0 \\ 
\varphi _{k}%
\end{array}%
\right),
 %\label{vectpropnuir}
\end{equation*}%
where $\psi_k$ is the unique solution of the following problem:

\begin{equation*} 
\left\{ 
	\begin{array}{ll}
\nu \psi_k''+ k^2\psi_k=q\varphi_k \: & \hbox{in }(0,\pi ), \\ 
	\noalign{\smallskip}
\psi_k(0)=\psi_k(\pi )=0. \\

\end{array}%
\right.
% \label{phikprob}
\end{equation*}%
then
  $$ \left(L^{\ast }-k^{2}I_{d}\right) \Phi _{1,k}^{\ast }=0 \quad \hbox{and} \quad \left( L^{\ast }-\nu k^{2}I_{d}\right) \Phi _{2,k}^{\ast }=0. $$
Moreover, an explicit expression of $\psi_{k}$ is given by :
\begin{equation} 
 \label{tpsiksolnuir}
\psi_{k}(x)=\psi'_{k}(0)\frac{\sqrt{\nu}}{k}\sin\left(\frac{kx}{\sqrt{\nu}}\right)+\frac{\sqrt{\nu}}{\nu k}\int_{0}^{x}\sin\left(\frac{k}{\sqrt{\nu}}(x-\xi)\right)q(\xi)\varphi_{k}(\xi)\,d\xi,
\end{equation}%
where 
\begin{equation} 
 \label{psi'k0nuir}
\psi'_{k}(0)=-\frac{\int_{0}^{\pi}q(s)\varphi_{k}(s)\sin\left(\frac{k}{\sqrt{\nu}}(\pi-s)\right)ds}{\nu \sin\left(\frac{k \pi}{\sqrt{\nu}}\right)}
\end{equation}%
\end{enumerate}

\end{prpstn}

\smallskip

 \begin{rmrk}
 \label{remark}
 Consider the following problem
 \begin{equation} 
\label{eqgene}
\left\{ 
	\begin{array}{ll}
y''+\lambda y =f(x)\quad \hbox{in  } (0,\pi),\\
\noalign{\smallskip}
y(0)= 0,\\
\noalign{\smallskip}
%y'(0)= \Theta
\end{array}%
\right.
\end{equation}%
where $\lambda \in \mathbb{R}^{\ast}_{+}$.
The general solution to (\ref{eqgene}) is given by
\begin{equation}
y(x)=a\frac{\sin(\sqrt{\lambda} x)}{\sqrt{\lambda}}+\frac{1}{\sqrt{\lambda}}\int_{0}^{x}\sin \left( \sqrt{\lambda}(x-s )\right) f(s)\,ds, \quad a \in \mathbb{R}^*,
\label{y}
\end{equation}
and
\begin{equation} 
\label{y'gene}
 y'(x)=a\cos(\sqrt{\lambda} x)+\int_{0}^{x}\cos \left( \sqrt{\lambda}(x-s )\right) f(s)\,ds,
\end{equation}
consequently
\begin{equation} 
\label{cnsygene}
y'(0)=a=-\frac{\int_{0}^{\pi}f(s)\sin\left(\sqrt{\lambda} (\pi-s)\right)ds}{\sin\left(\sqrt{\lambda}  \pi\right)}, \quad \hbox{if}\quad \sqrt{\lambda} \notin \mathbb{N^{\ast}}.
 \end{equation}%
 On the other hand, for all $\sqrt{\lambda} \notin \mathbb{N}^* $ :
 $$
 y(x)=\sum_{n\geq1}\frac{\int_{0}^{\pi }f( x )\varphi _{n}( x )\, dx}{\left(\lambda -n^2 \right) }\varphi _{n}( x ). 
$$
\end{rmrk}

\bigskip

\noindent
\textit {Proof of Proposition \ref{fam_nu_ir}}. 

Let us assume $\sqrt{\nu}  \notin \mathbb{Q}$. Let $ \lambda $ be an eigenvalue of $L^{\ast}$ and $y=(y_1,y_2)^{T}$ an associated eigenfunction. Thus $y$ is a solution of the following problem:

\begin{equation*} 
\left\{ 
	\begin{array}{ll}
	- y''_1=\lambda y_1 \: \hbox{in }\:(0,\pi ), \\ 
qy_1-\nu y''_2=\lambda y_2 \: \hbox{in }\: (0,\pi), \\ 
	\noalign{\smallskip}
y_1(0)=y_2(0)=0,\quad y_1(\pi)=y_2(\pi)=0.\\
\end{array}%
\right.  
\end{equation*}%
If $y_1\equiv 0$, then,  $\lambda=\nu k^2$ is an eigenvalue of $L^{\ast}$ and taking $y_2=\varphi_k$, we obtain $\Phi^{\ast}_{2,k}$ as associated eigenfunction of $L^{\ast}$. 
Now assume that $y_1\not\equiv 0$, then $\lambda=k^2$ and $y_1 = \varphi_k$ is a (normalized) solution to the first o.d.e. Inserting this expression in the second equation, we get for $y_2$ :
\begin{equation} 
\label{y2}
\left\{ 
	\begin{array}{ll}
y''_2+\frac{k^2}{\nu} y_2 =\frac{1}{\nu} q \varphi_k \: \hbox{in }\: (0,\pi), \\ 
	\noalign{\smallskip}	
y_2(0)=y_2(\pi)=0.\\
\end{array}%
\right.  
\end{equation}%
This proves that $\Phi _{1,k}^{\ast }$, is the second eigenfunction of $L^{\ast}$, associated to $k^2$ . Moreover (\ref{tpsiksolnuir}) (resp. (\ref{psi'k0nuir})) can be deduce from (\ref{y}) (resp. (\ref{cnsygene})).
$\hfill \square $

\begin{lmm}
\label{complete}
The sequence  $\mathcal{B^{\ast}}=\left\{ \Phi^{\ast} _{1,k},\Phi^{\ast} _{2,k} : k\in \mathbb{N}^{\ast} \right\}$  is complete in $L^{2}(0,\pi ;\mathbb{R}^2)$. 
\end{lmm}
\noindent
\textit {Proof of Lemma \ref{complete} }.\\
%\begin{enumerate}
Indeed, if $f=(f_1,f_2)$ is such that
$$
\langle f,\Phi_{\mu,k}^{\ast}\rangle=0, \quad \forall k \ge 1, \quad \forall \mu =1,2,
$$
then in particular
$$
\forall k \ge 1 \quad
\left\{  \begin{array}{l}
\langle f_2,\varphi_{k}\rangle=0 \\
\langle f_1,\varphi_{k}\rangle+\langle f_2,\psi_{k}\rangle =0.
\end{array} \right.
$$
This implies that $f_1=f_2=0$ (since $\{\varphi_{k} \}_{k\ge 1}$ is an orthonormal basis in $L^2(0,\pi)$ and proves the
completeness of $\mathcal{B^{\ast}}$. 
$\hfill \square $
  
 \begin{prpstn}
 	 There exists $\sqrt{\nu}\notin \mathbb{Q}^*_+$ and $q\in L^{\infty}(0,\pi) $ such that 
 	$\mathcal{B^{\ast}}=\left\{ \Phi^{\ast} _{1,k},\Phi^{\ast} _{2,k} : k\in \mathbb{N}^{\ast} \right\}$  is not a Riesz basis for $L^{2}(0,\pi ;\mathbb{R}^2)$.
 	\label{notriesz}
 \end{prpstn}
 \noindent
 To prove this result we will need the following two lemmas.
 %%%%%%%%%%%%%%%%%%%%%%%%%%%%lm12%%%%%%%%%%%%%
  \begin{lmm}
 For a sequence $\{f_k\}_{k \ge 1}$ in Hilbert space $(\mathcal{H},\langle ,\rangle)$ the following conditions are equivalent:
 \begin{enumerate}%[label=roman*)]
 \item $\{ f_k \}_{k \ge 1}$  is a Riesz basis for $\mathcal{H}$.
 \item $\{ f_k \}_{k \ge 1}$  is complete, and its Gram matrix $\left( \langle f_k,f_j \rangle \right)_{k,j \ge 1}$ defines a bounded, invertible operator on $l^2(\mathbb{N})$ the space of square summable scalar sequences.
 \end{enumerate}
 \label{lem}
 \end{lmm}
\noindent
\textit {Proof of Lemma \ref{lem} }. See for instance \cite{Ole}, Theorem 3.6.6, page 66.
 $\hfill \square $
 %%%%%%%%%%%%%%%%%%%%%%%%%%%%%%%%%%%%%%%%%%%%%%%%%
 \begin{lmm}
 \label{lem+}
% \begin{enumerate}%[label=roman*)]
% \item
% For any $\tau_0\in (0,\infty)$, there exist an irrational number $\nu > 0$ and a sequence of rational numbers $\{p_k,q_k \}_{k\ge 0}$ such that $p_k$ and $q_k$ are co-prime positive integers, the sequences $\{p_k\}_{k\ge 0}$ and $\{q_k \}_{k\ge 0}$ are strictly increasing and
%\begin{equation}
%\lim e^{\tau_0p_k^2}\left \vert \sqrt{\nu}-\frac{p_k}{q_k} \right \vert=1
%\label{lem62}
%\end{equation}
%  \item
  For any $\sigma \in (0,\infty)$, there exist an irrational number $\nu > 0$ and a sequence of rational numbers $\{k_p,j_p \}_{p\ge 0}$ such that $k_p$ and $j_p$ are co-prime positive integers, the sequences $\{k_p\}_{p\ge 0}$ and $\{j_p \}_{p\ge 0}$ are strictly increasing and
\begin{equation}
\lim e^{k_p^{2+\sigma}}\left \vert \sqrt{\nu}-\frac{k_p}{j_p} \right \vert=0.
\label{lem62+}
\end{equation}
In particular, we deduce the existence of a positive constant $C>0$
such that
\begin{equation}
\left \vert j_p\sqrt{\nu}-k_p \right \vert \le C k_p e^{-k_p^{2+\sigma}},\quad \forall p\ge 1.
\label{lem62+ en}
\end{equation}
% \end{enumerate}
 \end{lmm}
 \smallskip
 \noindent
\textit {Proof of Lemma \ref{lem+} }. See \cite{JFA14}, Lemma 6.22, page 47.
 $\hfill \square $
%%%%%%%%%%%%%%%%%%%%%%%%%%%%endlm12%%%%%%%%%%%%%%%%%%%%%%%%%%% 

\smallskip
\noindent
 We now have all the ingredients to prove Proposition \ref{notriesz}.
 
  \smallskip
  
 \noindent
\textit {Proof of Proposition \ref{notriesz} }.
The determinant of the Gram matrix associated to the normalized vectors of  $\mathcal{B^{\ast}}$ is equal to 
\begin{equation}
\begin{aligned}
 \det\, [G_{k_p,j_p}(\nu, q)]&= 1-\frac{\left \vert \langle \Phi_{1,k},\Phi_{2,j}\rangle \right \vert^2}{\Vert \Phi_{1,k}\Vert^2_{H^1_0}\Vert \Phi_{2,j}\Vert^2_{H^1_0}}\\
 &= 1-\frac{\left \vert \displaystyle\int _0^{\pi}\psi'_{k}(x)\varphi'_{j}(x)dx \right \vert^2}{j^2\left(k^2+\Vert \psi_k\Vert^2_{H^1_0}\right)} \\
 &= 1-\frac{j^2\left \vert \displaystyle\int _0^{\pi}q(x)\varphi_{k}(x)\varphi_{j}(x)dx \right \vert^2}{  \left( k^2+\displaystyle\underset{n\geq 1}\sum \frac{n^2 \left \vert\displaystyle\int _0^{\pi}q(x)\varphi_{k}(x)\varphi_{n}(x)dx \right \vert^2}{\left \vert k^2-\nu n^2 \right \vert^2}\right)\left \vert k^2-\nu j^2 \right \vert^2} \\
 &= 1-\frac{1}{  \left\lbrace 1+ U_{k,j}(\nu,q) \left[k^2+V_{k,j}(\nu,q)\right]\right\rbrace}, 
 \end{aligned}
 \label{det}
\end{equation}
 %we have used here the fact that : $\Vert \psi_k\Vert^2_{H^1_0}=\Vert \psi'_k\Vert^2_{L^2}=\sum_{n \ge 1} \left \vert \int _0^{\pi}\psi'_{k}(x)\varphi'_{n}(x)dx \right \vert^2$ 
where 
\begin{equation*}
U_{k,j}(\nu,q) =\frac{\left \vert k^2-\nu j^2 \right \vert^2}{j^2\left \vert \displaystyle\int _0^{\pi}q(x)\varphi_{k}(x)\varphi_{j}(x)dx \right \vert^2}, \quad V_{k,j}(\nu,q) =\displaystyle\underset{n\neq j}\sum \frac{n^2 \left \vert\displaystyle\int _0^{\pi}q(x)\varphi_{k}(x)\varphi_{n}(x)dx \right \vert^2}{\left \vert k^2-\nu n^2 \right \vert^2}.
%\label{lem62+}
\end{equation*}
Remark that, there exists a constant $C(\nu,\Vert q \Vert_{L^{\infty}})>0$ such that,% for all $\nu>0$ and $q\in L^{\infty}(0,\pi)$
\begin{equation*}
V_{k,j}(\nu,q) \le C(\nu,\Vert q \Vert_{L^{\infty}}),\quad  \forall k,j \ge 1.
%\label{V}
\end{equation*}
%\begin{crllr}
 Thanks to Lemma \ref{lem+}, we can extract a subsequence $(k_p+j_p)_{p\ge 0}$ of even numbers only or  odd numbers only, such that this two situations:
  \begin{enumerate}
  \item
  By choosing the subsequence of even numbers with $q(x)=\sin(x),\; x\in (0,\pi)$, we obtain 
  $$
  \left \vert\int _0^{\pi}q(x)\varphi_{k_p}(x)\varphi_{j_p}(x)dx \right \vert=\frac{2}{\pi}\frac{4k_pj_p}{\left [ (j_p -k_p)^2-1\right ]\left [ (j_p +k_p)^2-1\right ]}.
 $$
  Consequently
  $$
  k_p^2U_{k_p,j_p} \underset{p \longrightarrow +\infty}\longrightarrow 0,\quad \hbox{thus}\quad \det\, [G_{k_p,j_p}] \underset{p \longrightarrow +\infty}\longrightarrow 0.
  $$
   \item By choosing the subsequence of odd numbers with $q(x)=\sin(2x),\; x\in (0,\pi)$, we obtain 
  $$
  \left \vert\int _0^{\pi}q(x)\varphi_{k_p}(x)\varphi_{j_p}(x)dx \right \vert=\frac{2}{\pi}\frac{8k_pj_p}{\left \vert(j_p -k_p)^2-2\right \vert \left [ (j_p +k_p)^2-2\right ]}.
  $$
Consequently
  $$
  k_p^2U_{k_p,j_p} \underset{p \longrightarrow +\infty}\longrightarrow 0,\quad \hbox{thus}\quad \det\, [G_{k_p,j_p}] \underset{p \longrightarrow +\infty}\longrightarrow 0.
  $$
  \end{enumerate}
 % \label{coro}
%  \end{crllr}
%  \noindent
%\textit {Proof of Corollary \ref{coro} }
Indeed, take the first point (1) and let us fix $\sigma > 0$ and $\sqrt{\nu}>0$. We have

$$
\begin{aligned}
k_p\sqrt{ U_{k_p,j_p}}&= \frac{k_p\left \vert k_p^2-\nu j_p^2 \right \vert}{j_p\left \vert\displaystyle\int _0^{\pi}q(x)\varphi_{k_p}(x)\varphi_{j_p}(x)dx \right \vert}\\
 &=\frac{\pi}{2}\frac{k_p\left \vert k_p+\sqrt{\nu} j_p \right \vert\left [ (j_p -k_p)^2-1\right ]\left [ (j_p +k_p)^2-1\right ]}{2j^2_pk_p}\left \vert k_p-\sqrt{\nu} j_p\right \vert \\
 &\le \frac{\pi}{2} \frac{\left \vert k_p+\sqrt{\nu} j_p \right \vert\left [ (j_p -k_p)^2-1\right ]\left [ (j_p +k_p)^2-1\right ]}{2j^2_pk_p} k^2_p e^{-k_p^{2+\sigma}} \underset{p \longrightarrow +\infty}\longrightarrow 0. 
 \end{aligned}
$$
Thanks to formula (\ref{det}), 
$$\det\, [G_{k_p,j_p}] \underset{p \longrightarrow +\infty}\longrightarrow 0.$$
The point (2) can be similarly shown. 
 $\hfill \square $

\bigskip

Concerning the approximate controllability of system (\ref{boun.sg2x2}), it is well known that can be characterized in terms of a property of the solutions to (\ref{Padjoint}). More precisely,  system (\ref{boun.sg2x2}) is approximately controllable if and only if the following unique continuation property holds:

\smallskip

 \textit{"Let $\theta_0 \in H^{1}_0(0,\pi ;\mathbb{R}^{2})$  be given and let $\theta$  be the associated adjoint state. Then, if $B^{\ast }D\theta_x(0,t)=0$ on $(0,T)$, one has $\theta_0 \equiv 0 $  in $Q_T$"}.\\
Fattorini gave an interesting characterization of the approximate controllability under a general abstract framework. In his paper \cite{FAT66}, he proved that, under some reasonable assumptions, the only observation of the eigenfunctions completely characterizes the approximate controllability. Actually, this theorem has been proved for bounded observation operators but G. Olive (in \cite{GO}), give a generalization to the case of relatively bounded observation operators. We deduce that system (\ref{boun.sg2x2}) is approximately controllable at time
$T > 0$, if and only if for any $s\in \mathbb{C}$ and any $u \in \mathcal{D}(L^{\ast})=H^{2}(0,\pi ;\mathbb{R}^{2})\cap H^{1}_0(0,\pi ;\mathbb{R}^{2})$ we have 
\begin{equation} 
\label{FAT66}
\left.
\begin{aligned}
& L^{\ast}u = su, \\
&B^{\ast}D \partial_x u(0) = 0 
\end{aligned}
\right \} \Longrightarrow u \equiv 0  \quad \hbox{in}\; (0,\pi) .
\end{equation}
This previous relation (\ref{FAT66}) justifies the second point of the Theorem \ref{Thm} and will be used to prove the approximate controllability of system (\ref{boun.sg2x2}).

%%%%%%%%%%%%%%%%%%%%%%%%%%%%%%%%%%%%%%%%%%%%%%%%%%%%%
\subsection{Proof of Theorem \ref{Thm} : first point (a)}
%%%%%%%%%%%%%%%%%%

In this subsection, our objective is to prove that system  (\ref{boun.sg2x2}) is null controllable  at time $T$  if $T > \widetilde{T}_0\in [0,\infty)$, when $\sqrt{\nu}>1$. 
If $y$ is the solution of system (\ref{boun.sg2x2}) associated with $ y_{0}\in H^{-1} ( 0,\pi ;\mathbb{R}^{2}) $  and $u\in L^{2} ( 0,T ) $, then it can be checked that $y(T)=0$ in $(0,\pi) $ if and only if

\begin{equation}
\label{ipp}
\int_{0}^{T}  u(t)B^{\ast }D\theta_x(0,t) dt = - \langle y_0 , \theta (\cdot , 0) \rangle_{H^{-1},H_0^1} , \quad \forall
\theta _{0}\in H^{1}_0 ( 0,\pi ;\mathbb{R}^{2} ),
\end{equation}
where $\theta$ is the solution of the adjoint problem (\ref{Padjoint}) associated with $\theta_0$.
Taking 

$$
\theta _{0}:=\Phi _{i,k}^{\ast }, \; \forall k\geq 1,\; i=1,2,
$$ 
the corresponding solution to the adjoint problem (\ref{Padjoint}) is given by

\begin{equation*}
	\left\{ 
\begin{aligned}
&\theta^{1,k}(x, t)=e^{-k^2 ( T-t ) }\Phi _{1,k}^{\ast }(x), \\
& \theta^{2,k}(x, t)=e^{-\nu k^2 ( T-t ) }\Phi _{2,k}^{\ast }(x),  \; \forall k\geq 1.
 \end{aligned}%
\right.
\end{equation*}
Since the sequence $\mathcal{B^{\ast}}=\left\{ \Phi^{\ast} _{1,k},\Phi^{\ast} _{2,k} : k\in \mathbb{N}^{\ast} \right\}$  is complete in $L^{2}(0,\pi ;\mathbb{R}^2)$, the null controllability problem  for system (\ref{boun.sg2x2}) is equivalent to find $u\in L^2(0, T)$ such that:
%we reduce the null-controllability issue to the following problem
%of moments :
\begin{equation}
 \label{momInternnuir}
	\left\{ 
	\begin{aligned}
\int_0^T e^{-k^2t} u(T-t)\, dt &=-e^{-k^2T}\langle y_0,\psi_{1,k} \rangle_{H^{-1},H^1_0}
 ,  \\ 
\int_0^T e^{-\nu k^2t} u(T-t) \, dt &=-e^{-\nu k^2T}\langle y_0,\psi_{2,k} \rangle_{H^{-1},H^1_0}, \quad \forall k\geq 1,
\end{aligned}%
\right.  
\end{equation}
where $\psi_{1,k}$ and $\psi_{2,k}$ are defined in (\ref{psik}) . We are now going to give  some results that  will be crucial, to solve (\ref{momInternnuir}). One has:

\begin{prpstn}
\label{epp}
Let $\sqrt{\nu}>1$. Let us define 
$$
\begin{array}{ccc}
I : &\mathbb{N}^{\ast}&\longrightarrow \mathbb{N}^{\ast} \\
&k& \longmapsto i_k
\end{array}
$$
where, for any $k \in \mathbb{N}^*$, $i_k$ is the nearest integer to $\sqrt{\nu}k$  i.e 
$$\left \vert \sqrt{\nu}k -i_k   \right \vert<\frac{1}{2}.$$
Thus for any $k \in \mathbb{N}^*$,
$$\vert \sqrt{\nu}k-i \vert >\frac{1}{2},\quad  \forall \,i \in \mathbb{N}^*,\, i\neq i_k.$$
Then
\begin{enumerate}
\item 
The function $I$ is injective.
 \item 
$
 \widehat{I}=\mathbb{N}^*\setminus I(\mathbb{N}^*)=\{ \widehat{i}_k: k\ge 1\} 
$, is a infinite set, where the elements of $\widehat{I}$ are classified in ascending order.
\end{enumerate}
\end{prpstn}

\smallskip

\noindent
\textit {{Proof of Proposition \ref{epp}}}.
\begin{enumerate}
\item
 Let us assume that $I(k_1)=I(k_2)$, where $k_1$, $k_2$ $\in$ $\mathbb{N}^{\ast}$. We have $i_{k_1}=i_{k_2}$ with
$\left \vert \sqrt{\nu}k_1 -i_{k_1}   \right \vert<\frac{1}{2}$ and $\left \vert \sqrt{\nu}k_2 -i_{k_2}   \right \vert<\frac{1}{2}$. This leads to $$-1+(i_{k_1}-i_{k_2} )<\sqrt{\nu}(k_1-k_2)< 1+(i_{k_1}-i_{k_2} ),$$ and thus $\left \vert k_1-k_2 \right \vert < \frac{1}{\sqrt{\nu}}<1$, which implies $k_1=k_2$.

\item
Assume  $\sqrt{\nu}>1$, $\sqrt{\nu}\notin \mathbb{Q^*_+} $.
\begin{itemize}
	\item If $\sqrt{\nu} > 2$ then,  $\forall n\in \mathbb{N^*}$, $i_{n+1}-i_n>1$. %i.e $\exists k_0 \ge 1 $ such that $i_n<\widehat{i}_{k_0}<i_{n+1}$.
	\item  Suppose $1 <\sqrt{\nu} < 2$.	There exists a sequence of integers $(n_k)_{k\in \mathbb{N}^*}$ strictly increasing,  such that	$i_{n_{k}+1}-i_{n_k}>1$. Actually, let us take for instance
	$$ n_k \in \left]\frac{2k+1}{2(\sqrt{\nu}-1)}-1, \frac{2k+1}{2(\sqrt{\nu}-1)}\right [,\quad k\ge 1 .$$
	We deduce 
	$$
	-\frac{1}{2}<\frac{1}{2}-\sqrt{\nu}+1<n_k\sqrt{\nu}-(n_k+k)< \frac{1}{2},\quad i.e \quad i_{n_k}=n_k+k.
	$$
	Moreover
	 $$
	 \begin{array}{ccc}
	 i_{n_{k}+1}-i_{n_k}= i_{n_{k}+1}-n_k-k &>&\sqrt{\nu}(n_k+1)-\frac{1}{2}-n_k-k \\
	 &=& n_k\sqrt{\nu}-(n_k+k)+\sqrt{\nu}-\frac{1}{2}\\
	 &>& \frac{1}{2}-\sqrt{\nu}+1+\sqrt{\nu}-\frac{1}{2}=1.
	  \end{array}
	 $$
		
	\end{itemize}
	\end{enumerate}
$\hfill \square $

%%%%%%%%%%%%%%%%%%%%%%%%%%%%%%%%%%%%%%%%%%%%%%%%%
\subsubsection  {Positive null controllability result}
%%%%%%%%%%%%%%%%%%
Thanks to Proposition \ref{epp} we can reformulate (\ref{momInternnuir}). We say that the null controllability property at time $T$ for system (\ref{boun.sg2x2}) is equivalent to find $u\in L^2(0, T)$ such that:

\begin{equation}
 \label{momInternnuirbis}
	\left\{ 	
	\begin{aligned}
	 \int_0^T e^{-\widehat{i}_k^2t} u(T-t) \, dt &=-e^{-\widehat{i}_k^2T}\langle y_0,\psi_{1,\widehat{i}_k} \rangle_{-1,1},\\
\int_0^T e^{-i_k^2t} u(T-t)\, dt &=-e^{-i_k^2T}\langle y_0,\psi_{1,i_k} \rangle_{-1,1}
 ,  \\ 
\int_0^T e^{-\nu k^2t} u(T-t) \, dt &=-e^{-\nu k^2T}\langle y_0,\psi_{2,k} \rangle_{-1,1},
\end{aligned}
\quad \forall k\ge 1.
\right.  
\end{equation}
We can now state our following main result.

\begin{prpstn}\label{control}

Let us introduce the (closed) space $E_T \subset L^2(0, T)$ given by
$$
E_T=\overline{span\{e^{-k^2t},e^{-\nu k^2t}: k\ge1\}}^{L^2(0,T)}.
$$
 Then
\begin{enumerate}
\item
There exists a family $\{q_k \}_{k \ge 1} \subset E_T$  such that 

\begin{equation}
 \label{biorth}
 \left\{ 
\begin{aligned}
&	 \int_0^T e^{-\widehat{i}_k^2t} q_{k}(t) \, dt=-e^{-\widehat{i}_k^2T}\langle y_0,\psi_{1,\widehat{i}_k} \rangle_{-1,1},\\
&\int_0^T e^{-i_k^2t} q_{k}(t)\, dt=-e^{-i_k^2T}\langle y_0,\psi_{1,i_k} \rangle_{-1,1}
 ,  \\ 
&\int_0^T e^{-\nu k^2t} q_{k}(t) \, dt=-e^{-\nu k^2T}\langle y_0,\psi_{2,k} \rangle_{-1,1},\quad k\ge 1,
\\
&     \int_0^T  e^{-\widehat{i}_k^2t}q_{j}(t)dt=\int_0^T  e^{-i_k^2t}q_{j}(t)dt=\int_0^T  e^{-\nu k^2 t}q_{j}(t)dt=0, \quad k\neq j.
    \end{aligned}
 \right.
\end{equation}

\item If $T>\widetilde{T}_0=\max \left \lbrace T_1,\widetilde{T}_2 \right \rbrace$ (see (\ref{T0tilde})) then we infer that an explicit solution $u$  of moment problem (\ref{momInternnuirbis}) is given by

\begin{equation*}
%\label{unuir}
u(t)=\sum_{k\ge 1}q_k(T-t).
\end{equation*}
\end{enumerate}
\end{prpstn}

\noindent
\textit {{Proof of Proposition \ref{control}}}.
Let us start by recalling classical properties of the Laplace transform (see
for instance \cite{SHA}, pp.19-20). Let  $H^2(\mathbb{C}_+)$ the space of holomorphic functions $\Phi$ on $\mathbb{C}_+=\{ z \in \mathbb{C}, \Re(z)>0\}$ such that
$$
\sup_{\sigma>0}\Vert \Phi(\sigma +i \bullet)\Vert_{L^2(\mathbb{R}; \mathbb{C})}<\infty, 
$$
endowed with the norm

$$
\Vert \Phi\Vert^2_{H^2(\mathbb{C}_+)}=\sup_{\sigma>0}\Vert \Phi(\sigma +i \bullet)\Vert^2_{L^2(\mathbb{R}; \mathbb{C})}=\int_{\mathbb{R}}|\Phi(i \tau)|^2 d\tau .
$$
Then the Laplace transform
$$
\begin{array}{ccc}
L  :L^2(0,+\infty; \mathbb{C})&\longrightarrow& H^2(\mathbb{C}_+)\\
f  &\longmapsto& \left( \Phi\ : \lambda \in \mathbb{C^+}\longmapsto \displaystyle  \int_{\mathbb{R}}e^{-\lambda t}f(t)dt \right).
\end{array}
$$
is an isomorphism. In the sequel, We shall construct for each $k$, a function $J_k \in H^2(\mathbb{C}_+)$ satisfying appropriate properties, and take advantage of the isomorphism property of the Laplace transform to conclude. 
Let us fix $k\ge 1$, we define $J_k$, for $\lambda \in \mathbb{C}_+$, as
\begin{equation*}  %\label{jk}
\displaystyle
J_k(\lambda)= \left\{ \alpha_k(\lambda-i_k^2)(\lambda-\widehat{i}_k^2)+ \beta_k(\lambda-\nu k^2)(\lambda-\widehat{i}_k^2)+\gamma_k(\lambda-\nu k^2)(\lambda-i_k^2)\right\} L_k(\lambda),
\end{equation*}
where $\alpha_{k}, \beta_k$ and $\gamma_k$ are constants to be determined and $L_k$ is the
 following Blaschke-type product 
\begin{equation}  \label{Lk}
L_k(\lambda)=\frac{1}{(1+\lambda)^3}\displaystyle{\prod_{ j\ge1, j\neq k  } } \left( \frac{\lambda-\nu j^2}{\lambda+\nu j^2} \right) \left(\frac{\lambda-i_j^2}{\lambda+i_j^2} \right)\left(\frac{\lambda-\widehat{i}_j^2}{\lambda+\widehat{i}_j^2}\right).
\end{equation}
Notice that for any $k\ge 1$,   
\begin{equation*}
\left \vert L_k(i \tau) \right \vert^2 =\frac{1	}{(1+\tau^2)^3}, \; \forall\; \tau \in \mathbb{R}.
\end{equation*}
This implies $J_k \in H^2(\mathbb{C}_+) $, moreover
\begin{equation*}  %\label{jk+}
J_k(\nu j^2)=J_k(i_j^2)=J_k(\widehat{i}_j^2)=0,\quad \forall j \ge 1,\, j\neq k.
\end{equation*}
So, using that the Laplace transform is a isomorphism from $L^2(0,\infty;\mathbb{C})$ into $H^2(\mathbb{C}_+)$, we infer the existence of a nontrivial function $\tilde{q}_k \in L^2(0,\infty;\mathbb{C}) $ such that
\begin{equation*} 
%\label{iso}
J_k(\lambda)= \int_{0}^{+\infty}e^{-\lambda t}\widetilde{q}_k(t)dt,\quad \forall \lambda \in \mathbb{C}_+.
\end{equation*}
Now, we can choose $\alpha_k,\, \beta_k$ and $\gamma_k$ such that 

\begin{equation*}
%\label{alphakbetak}
\left\{ 
	\begin{aligned}
&J_k(\nu k^2)=-e^{-\nu k^2T}\langle y^0,\psi_{2,k} \rangle_{-1,1}, \\ 
&J_k(i_k^2)=-e^{-i_k^2T}\langle y^0,\psi_{1,i_k} \rangle_{-1,1}, \\ 
&J_k(\widehat{i}_k^2)=-e^{-\widehat{i}_k^2T}\langle y^0,\psi_{1,\widehat{i}_k} \rangle_{-1,1},
\end{aligned}%
\quad \forall k\geq 1.
\right. 
\end{equation*}
We obtain 
\begin{equation}
 \label{momIntern+}
	\left\{ 
	\begin{aligned}
&\alpha_k =-\frac{e^{-\nu k^2T}}{(\nu k^2-i_k^2)(\nu k^2-\widehat{i}_k^2) L_k(\nu k^2) }\langle y_0,\psi_{2,k} \rangle_{-1,1},  \\ 
&\beta_k =-\frac{e^{-i_k^2T}}{(i_k^2- \nu k^2)(i_k^2- \widehat{i}_k^2) L_k(i_k^2)}\langle y_0,\psi_{1,i_k}\rangle_{-1,1},\\
&\gamma_k =-\frac{e^{-\widehat{i}_k^2T}}{(\widehat{i}_k^2- \nu k^2)(\widehat{i}_k^2- i_k^2) L_k(\widehat{i}_k^2)}\langle y_0,\psi_{1,\widehat{i}_k}\rangle_{-1,1}, 
\end{aligned}%
\quad \forall k\geq 1.
\right.  
\end{equation}
The Parseval equality gives

$$
\begin{aligned}
\Vert \widetilde{q}_k \Vert^2_{L^2(0,\infty;\mathbb{C})}&= \displaystyle\int_{-\infty}^{+\infty}|J_k(i\tau)|^2 d \tau,\quad \quad \quad \\
%\displaystyle  &=& \displaystyle \int_{-\infty}^{+\infty}\left \vert (i \tau- \widehat{i}_k^2)\left \{i\tau(\alpha_k+\beta_k)-(\alpha_k i_k^2 +\beta_k\nu k^2)\right\} +\gamma_k(i \tau- \nu k^2)(i \tau- i_k^2) \right \vert^2|L_k(i\tau)|^2 d \tau,\\
&= \int_{-\infty}^{+\infty}\frac{\left \vert (i \tau- \widehat{i}_k^2)\left \{(\alpha_k+\beta_k)(i \tau-i_k^2)+ \beta_k (i_k^2-\nu k^2)\right\} +\gamma_k(i \tau- \nu k^2)(i \tau- i_k^2) \right \vert^2}{\vert 1+i\tau\vert^6} d \tau,\\
&=  \int_{-\infty}^{+\infty}\frac{\left \vert (\alpha_k+\beta_k)(i \tau-i_k^2)(i \tau- \widehat{i}_k^2)+ \beta_k (i_k^2-\nu k^2)(i \tau- \widehat{i}_k^2) +\gamma_k(i \tau- \nu k^2)(i \tau- i_k^2) \right \vert^2}{\vert 1+\tau^2\vert^3}d \tau,\\
&\le  C\nu^2 k^4 i_k^4\widehat{i}_k^4\left \{ (\alpha_k+\beta_k)^2+ \beta_k^2 (i_k^2-\nu k^2)^2+  \gamma_k^2  \right\}.
\end{aligned}
$$
In the sequel, we will show that there exists $q_k\in L^2(0,T)$ such that 

\begin{equation}
\label{normqk}
\Vert q_k \Vert^2_{L^2(0,T)}
\le  C \nu^2 k^4 i_k^4\widehat{i}_k^4\left \{ (\alpha_k+\beta_k)^2+ \beta_k^2 (i_k^2-\nu k^2)^2+  \gamma_k^2  \right\}.
\end{equation}
Assume that a countable family $\Lambda$ satisfies the \textit{weak gap condition}\footnote{This notion was introduced in   \cite{BBM18}}. On one hand, let us consider the closed space $A(\Lambda,T)\subset L^2(0,T; \mathbb{C})$ given by
    $$
    A(\Lambda,T)=\overline{\hbox{Span}\{ e^{-\lambda t}; \, \lambda \in \Lambda\}}^{L^2(0,T;\mathbb{C})}.
    $$
    For any $T\in (0,+\infty)$, the restriction operator
    \begin{equation}
    \label{restricop}
    R_{\Lambda,T}: \phi \in A(\Lambda,\infty) \rightarrow R_{\Lambda,T} \phi=\phi_{|(0,T)} \in A(\Lambda,T)
    \end{equation}
    is an isomorphism. Moreover there exists a constant $C_{T}$ such that
    \begin{equation}
    \label{restrictionoperator}
    \Vert R^{-1}_{\Lambda,T} \Vert \le C_{T}.
    \end{equation}
Indeed, the fact that $R_{\Lambda,T}$ is an isomorphism is proved in \cite{SCH}. The proof of the bound (\ref{restrictionoperator}) can be done by contradiction (see \cite{JFA14} Lemma 4.2 and \cite{BBM18} Proposition 2.4 ).   \\
On the other hand, for any $T > 0$, there exists a constant $C_{T}$ such that for any $\widetilde{q} \in L^2(0,+\infty;\mathbb{C})$ there exists $q \in L^2(0,T;\mathbb{C})$ satisfying
$$
\int_0^T q(t)e^{-\lambda t}dt = \int_0^{+\infty} \widetilde{q}(t)e^{-\lambda t}dt,\quad \forall \lambda \in \Lambda,
$$
and 
$$
\left \Vert q\right \Vert_{L^2(0,T,\mathbb{C})}\le C_{T}\left \Vert \widetilde{q}\right \Vert_{L^2(0,+\infty,\mathbb{C})}.
$$
Indeed, from \cite{JMPA11} Corollary 4.3, it comes that $A(\Lambda,\infty)$ is a proper subspace of $L^2(0,+\infty;\mathbb{C})$. Let $\Pi_{\Lambda}$ the associated orthogonal projection and  $ \widetilde{q} \in L^2(0,+\infty;\mathbb{C})$.  Then, by construction, we have
 
\begin{equation}
 \label{const_q}   
\int_0^{+\infty} \Pi_{\Lambda} \widetilde{q}(t)e^{-\lambda t}dt=\int_0^{+\infty}  \widetilde{q}(t)e^{-\lambda t}dt, \quad \forall \lambda \in \Lambda.
\end{equation}
Since the restriction operator $R_{\Lambda,T}$ defined in (\ref{restricop}) is an isomorphism, choose  $q:=(R_{\Lambda,T}^{-1})^*\Pi_{\Lambda}\widetilde{q}$, thus, there exists $C_{T}$ such that

$$
\left \Vert q\right \Vert_{L^2(0,T,\mathbb{C})}\le C_{T}\left \Vert \widetilde{q}\right \Vert_{L^2(0,+\infty,\mathbb{C})}.
$$
Using  (\ref{const_q}), for  every $\lambda \in \Lambda$
$$
\begin{aligned}
\int_0^{T} q(t)e^{-\lambda t}dt&=\langle (R_{\Lambda,T}^{-1})^*\Pi_{\Lambda}\widetilde{q}(\cdot), e^{-\lambda \cdot} \rangle _{L^2(0,T)}, \\
&=\langle \Pi_{\Lambda}\widetilde{q}(\cdot), R_{\Lambda,T}^{-1}R_{\Lambda,T}e^{-\lambda \cdot} \rangle _{L^2(0,+\infty)},\\
&= \int_0^{T} \widetilde{q}_k(t)e^{-\lambda t}dt.
\end{aligned}
$$
This proves inequality (\ref{normqk}).\\
Now, thanks to equalities (\ref{momIntern+}), we obtain for all $ k\ge 1 $:
\begin{equation*}
%\label{alpha+beta}
\left\{ 
	\begin{aligned}
&(\alpha_k+\beta_k)^2 = \frac{1}{(\nu k^2-i_k^2)^2} \left \vert \left\langle y^0,\frac{e^{-\nu k^2T}\psi_{2,k} }{(\nu k^2-\widehat{i}_k^2)  L_k(\nu k^2)  }-\frac{e^{-i_k^2T}\psi_{1,i_k} }{(i_k^2-\widehat{i}_k^2) L_k(i_k^2)  },\right\rangle_{-1,1}  \right \vert^2,
\\
&\beta_k^2 (i_k^2-\nu k^2)^2+  \gamma_k^2=\frac{1}{(\widehat{i}_k^2-i_k^2)^2}\left( \frac{e^{-2i_k^2T}  \left \vert \left\langle y^0,\psi_{1,i_k} \right\rangle_{-1,1}  \right \vert^2}{\left\vert L_k(i_k^2)\right\vert^2  }+\frac{e^{-2\widehat{i}_k^2T}\left \vert \left\langle y^0,\psi_{1,\widehat{i}_k}\right\rangle_{-1,1}  \right \vert^2}{(\widehat{i}_k^2- \nu k^2)^2 \left \vert L_k(\widehat{i}_k^2)\right\vert^2} \right).
\end{aligned}%
\right. 
\end{equation*}
Thus the control $u$ defines an element of $L^2(0,T)$ if
\begin{equation}
\label{serie1}
U=\displaystyle{\sum_{k\ge1}}U_k=\displaystyle{\sum_{k\ge1}} \frac{(k i_k\widehat{i}_k)^4}{(\nu k^2-i_k^2)^2} \left \vert \left\langle y^0,\frac{e^{-\nu k^2T}\psi_{2,k} }{(\nu k^2-\widehat{i}_k^2)  L_k(\nu k^2)  }-\frac{e^{-i_k^2T}\psi_{1,i_k} }{(i_k^2-\widehat{i}_k^2)  L_k(i_k^2)  }\right\rangle_{-1,1}  \right \vert^2<+\infty,
\end{equation}
and
\begin{equation}
\label{serie2}
V=\displaystyle{\sum_{k\ge1}}V_k=\displaystyle{\sum_{k\ge1}} \frac{(k i_k\widehat{i}_k)^4}{(\widehat{i}_k^2-i_k^2)^2}\left( \frac{e^{-2i_k^2T}  \left \vert \left\langle y^0,\psi_{1,i_k} \right\rangle_{-1,1}  \right \vert^2}{\left\vert L_k(i_k^2)\right\vert^2  }+\frac{e^{-2\widehat{i}_k^2T}\left \vert \left\langle y^0,\psi_{1,\widehat{i}_k}\right\rangle_{-1,1}  \right \vert^2}{(\widehat{i}_k^2- \nu k^2)^2 \left \vert L_k(\widehat{i}_k^2)\right\vert^2} \right)<+\infty.\\
\end{equation}
Let us give first some properties of $L_k$ (see (\ref{Lk})).

\begin{lmm}
\label{propoLk}
Assume  $\sqrt{\nu}>1$, $\sqrt{\nu}\notin \mathbb{Q^*_+} $.
\begin{enumerate}
\item
For every $\varepsilon > 0$ there exists positive constants $C_1(\varepsilon )$, $C_2(\varepsilon )$, and $C_3(\varepsilon )$, independent of $k$, such that
\begin{equation*}
%\label{estimLk2}
%\left\{ 
%	\begin{array}{l}
|L_k(i_k^2)|\ge C_1(\varepsilon )e^{-\varepsilon i_k^2},\:
|L_k(\widehat{i}_k^2)|\ge C_2(\varepsilon )e^{-\varepsilon \widehat{i}_k^2},\hbox{  and  }
|L_k(\nu k^2)|\ge C_3(\varepsilon )e^{-\varepsilon \nu k^2}, \:\forall k\ge 1.
%\right. 
\end{equation*}
\item
 There exists a constant $C> 0$, independent of $k$, such that
\begin{equation*}
%\label{estimL'k}
|L'_k(\lambda)|\le C, \quad \forall \lambda \ge 0.
\end{equation*}
\end{enumerate}
\end{lmm}

\noindent
\textit {Proof of Lemma \ref{propoLk}}. 
\begin{enumerate}

\item
Let us fix $k \ge 1 $ and work with $L_k(\nu k^2)$. We have
\begin{equation}
\label{lknuk2}
|L_k(\nu k^2)|=\frac{1}{(1+\nu k^2)^3}\displaystyle{\prod_{ \begin{array}{ccc} j\ge1\\ j\neq k  \end{array}} }\frac{ \vert 1-\frac{\nu k^2}{\nu j^2} \vert}{(1+\frac{\nu k^2}{ \nu j^2})}\frac{\vert 1-\frac{ \nu k^2}{ i_j^2}\vert}{(1+\frac{ \nu k^2}{ i_j^2})}\frac{\vert 1-\frac{ \nu k^2}{ \widehat{i}_j^2}\vert}{(1+\frac{\nu k^2}{ \widehat{i}_j^2})}.
\end{equation}

$\bullet$ Let us fix $\varepsilon>0$,  there exists $N(\varepsilon)$ such that 

$$
{\sum_{ j\ge N(\varepsilon)} \left(\frac{1}{\nu j^2}+\frac{1}{i_j^2}+\frac{1}{\widehat{i}_j^2}\right)} \le \varepsilon.
$$
Using the inequality $1+x\le e^x$, $x\in \mathbb{R}$, we can estimate the denominator of (\ref{lknuk2}) as follows:
\begin{equation*}
\begin{aligned}
&(1+\nu k^2)^3\displaystyle\prod_{  j\ge1, j\neq k   }\left(1+\frac{\nu k^2}{\nu j^2}\right)\left(1+\frac{ \nu k^2}{ i_j^2}\right)\left(1+\frac{ \nu k^2}{ \widehat{i}_j^2}\right) \\
&\le(1+\nu k^2)^3 \displaystyle{\prod_{ j=1}^{N(\varepsilon)-1}} \left(1+\nu k^2\right)^3\displaystyle{\prod_{ j\ge N(\varepsilon)}} \left(1+\frac{\nu k^2}{\nu j^2}\right)\left(1+\frac{ \nu k^2}{ i_j^2}\right)\left(1+\frac{ \nu k^2}{ \widehat{i}_j^2}\right)\\
&\le (1+\nu k^2)^3 \displaystyle{\prod_{ j=1}^{N(\varepsilon)-1}} \left(1+\nu k^2\right)^3\exp\left\{\nu k^2\displaystyle{\sum_{ j\ge N(\varepsilon)} \left(\frac{1}{\nu j^2}+\frac{1}{i_j^2}+\frac{1}{\widehat{i}_j^2}\right)} \right\}\\
&\le (1+\nu k^2)^{3N(\varepsilon)} e^{\varepsilon \nu k^2}\le   C(\varepsilon) e^{ \varepsilon \nu k^2 } ,\quad \forall k\in \mathbb{N}^*,
\end{aligned}
\end{equation*}

for a positive constant $C(\varepsilon)$.

$\bullet$ Let us now work on the numerator of (\ref{lknuk2}). 

Let us recall that (see Proposition \ref{epp} ),  for all $k,j\ge 1$ with $j\neq k$, we have

$$
\vert \sqrt{\nu}j-\sqrt{\nu}k\vert \ge \sqrt{\nu}, \quad \vert i_j-\sqrt{\nu}k \vert > \frac{1}{2}\quad \hbox{and} \quad  \vert \widehat{i}_j-\sqrt{\nu}k \vert > \frac{1}{2},
$$
 we deduce that, the conditions (1.8) et (1.9) of Theorem 1.1 in \cite{FATRUS}  are satisfied. Thus for every $\varepsilon > 0$ there exists positive constants $C_1(\varepsilon )$, $C_2(\varepsilon )$ and $C_3(\varepsilon )$ such that
$$
%\left\{
%\begin{aligned}
\displaystyle{\prod_{ j\ge1,j\neq k  } }\left\vert 1-\frac{\nu k^2}{\nu j^2}\right\vert \ge  C_1(\varepsilon )e^{-\varepsilon \nu k^2},\,
\displaystyle{\prod_{ j\ge 1,j\neq k  } }\left\vert 1-\frac{\nu k^2}{i_j^2}\right\vert \ge  C_2(\varepsilon )e^{-\varepsilon \nu k^2}\hbox{ and } \displaystyle{\prod_{ j\ge 1,j\neq k  } }\left\vert 1-\frac{\nu k^2}{i_j^2}\right\vert \ge  C_3(\varepsilon )e^{-\varepsilon \nu k^2}.
%\end{aligned}
%\right.
$$
\item
For all $\lambda \ge 0$, one has
$$
%\begin{aligned}
\frac{L'_k(\lambda)}{L_k(\lambda)}=\frac{d}{d\lambda}\log L_k(\lambda)=  \frac{-3}{1+\lambda}+\displaystyle{\sum_{ j\ge 1,j\neq k}}\left (\frac{2\widehat{i}_j^2}{(\lambda+\widehat{i}_j^2)^2}\frac{\lambda+\widehat{i}_j^2}{\lambda-\widehat{i}_j^2} +\frac{2\nu j^2}{(\lambda+\nu j^2)^2}\frac{\lambda+\nu j^2}{\lambda-\nu j^2}+\frac{2i_j^2}{(\lambda+i_j^2)^2}\frac{\lambda+i_j^2}{\lambda-i_j^2}\right),
%\end{aligned}
$$
consequently
$$
%\begin{array}{ccc}
\left \vert L'_k(\lambda)\right \vert
\le 3+2\displaystyle{\sum^{+\infty}_{j=1}}\left( \frac{1}{\widehat{i}_j^2}+ \frac{1}{\nu j^2}+\frac{1}{i_j^2}\right) <+\infty.
%\end{array}
$$
\end{enumerate}
$\hfill \square $

Let us now show that $V$ (see (\ref{serie2})) is a convergent series,  when $T>T_1$ (see (\ref{T0tilde})). One has

  $$   
\begin{aligned}
 V
&\le C \displaystyle{\sum_{k\ge1}} k^4 i_k^4\widehat{i}_k^4\left( \frac{e^{-2i_k^2T}  \left \vert \left\langle y_0,\psi_{1,i_k} \right\rangle_{-1,1}  \right \vert^2}{\left\vert L_k(i_k^2)\right\vert^2  }+\frac{e^{-2\widehat{i}_k^2T}\left \vert \left\langle y_0,\psi_{1,\widehat{i}_k}\right\rangle_{-1,1}  \right \vert^2}{\left\vert L_k(\widehat{i}_k^2)\right\vert^2} \right)\\
 &\le C(\varepsilon)\Vert y_0 \Vert^2_{-1} \displaystyle{\sum_{k\ge1}} k^4 i_k^4\widehat{i}_k^4\left( e^{-2i_k^2(T-\varepsilon)}  \left \Vert \psi_{1,i_k} \right \Vert^2+e^{-2\widehat{i}_k^2(T-\varepsilon)}\left \Vert \psi_{1,\widehat{i}_k}\right \Vert^2\right)\\
& \le C(\varepsilon)\Vert y_0 \Vert^2_{-1} \displaystyle{\sum_{k\ge1}}  k^4 i_k^4\widehat{i}_k^4e^{-2k^2(T-\varepsilon)}  \left \Vert \psi_{1,k} \right \Vert^2\\
& = C(\varepsilon)\Vert y_0 \Vert^2_{H^-1} \displaystyle{\sum_{k\ge1}}  k^4 i_k^4\widehat{i}_k^4e^{-2k^2(T-\frac{\log  \left \Vert \psi_{1,k} \right \Vert}{k^2}-\varepsilon)} \\
 & = C(\varepsilon)\Vert y_0 \Vert^2_{H^-1} \displaystyle{\sum_{k\ge1}}  k^4 i_k^4\widehat{i}_k^4e^{-2k^2(T-T_1-2\varepsilon)} \\
\end{aligned}
$$
If we take, $\varepsilon\in (0,\frac{T-T_1}{2})$, we deduce that $V$  converges if $T>T_1$.
Let us see now that $U$ (see  (\ref{serie1})) is a convergent series,  when $T>\widetilde{T}_2$ (see (\ref{T0tilde})). One has, for all $k\ge 1$s

  $$   
\begin{aligned}
 U_k &= \frac{k^4i_k^4\widehat{i}_k^4}{(\nu k^2-i_k^2)^2} \left \vert \left\langle y_0,\frac{e^{-i_k^2T}\psi_{1,i_k} }{(i_k^2-\widehat{i}_k^2)  L_k(i_k^2)  }-\frac{e^{-\nu k^2T}\psi_{2,k} }{(\nu k^2-\widehat{i}_k^2)  L_k(\nu k^2)  }\right\rangle\right \vert^2, \\
 &\le   \frac{k^4i_k^4\widehat{i}_k^4 \Vert y_0 \Vert^2_{-1}}{\left ( \nu k^2-i_k^2\right)^2}\left\Vert \frac{e^{-i_k^2T}\left(\psi_{1,i_k}-\psi_{2,k} \right)}{(i_k^2-\widehat{i}_k^2)  L_k(i_k^2)}+
\left (\frac{e^{-i_k^2T} }{(i_k^2-\widehat{i}_k^2)  L_k(i_k^2) }-\frac{e^{-\nu k^2T}}{(\nu k^2-\widehat{i}_k^2)  L_k(\nu k^2) }\right) \psi_{2,k}\right\Vert^2, \\
%&\le& k^4 i_k^4\widehat{i}_k^4 C(\varepsilon)\Vert y_0 \Vert^2_{-1}  \left\lbrace\frac{e^{-2i_ k^2(T-\varepsilon)} \left\Vert \psi_{1,i_k} -\psi_{2,k}\right\Vert^2}{\left ( \nu k^2-i_k^2\right)^2} 
% +  e^{-2\lambda_k^*T+\varepsilon i^2_k}  +e^{-2\nu k^2(T-\varepsilon)} \left( e^{\varepsilon i^2_k}\left \vert L'_k( \lambda_k^*)\right \vert^2 +1\right)  \right\rbrace,\\
 %\\
&\le   k^4 i_k^4\widehat{i}_k^4 C(\varepsilon)\Vert y_0 \Vert^2_{-1}  \left( \exp\left\{-2i_ k^2\left[ T-\frac{\log \frac{\left \Vert \psi_{2,k} -\psi_{1,i_k} \right \Vert}{\left\vert  \nu k^2-i_k^2 \right\vert }}{ i_k^2} -\varepsilon\right]\right\}
+    e^{-2\lambda_k^*(T-\frac{i_k^2}{\lambda_k^*}\varepsilon)}  +
e^{-2\nu k^2(T-\varepsilon-\frac{i_k^2}{\nu k^2}\varepsilon)} 
 \right),\\ 
&\le   k^4i_k^4\widehat{i}_k^4 C(\varepsilon')\Vert y_0 \Vert^2_{-1}  \left( e^{-2i_ k^2\left[ T-\widetilde{T}_2 -\varepsilon'\right]}
+    e^{-2\lambda_k^*(T-2\varepsilon')}  +e^{-2\nu k^2(T-\varepsilon')} 
 \right), 
\end{aligned}
$$
where  $\lambda^{\ast}_{k}=(1-\theta_k)\nu k^2+\theta_k i_k^2$, with $\theta_k\in ]0,1[$. If we take, $\varepsilon'\in (0,\frac{T-\widetilde{T}_2}{2})$, we deduce that $U$   converge if $T>\widetilde{T}_2$. This gives the proof of null-controllability of system (\ref{boun.sg2x2}) if $T>\widetilde{T}_0=\max \left \lbrace T_1,\widetilde{T}_2 \right \rbrace$.
$\hfill \square $

%%%%%%%%%%%%%%%%
\subsubsection  {The negative null controllability result}
%%%%%%%%%%%%%%
Let us prove that if $0<T<\widetilde{T}_0$, then system (\ref{boun.sg2x2}) is not null controllable at time $T$. We argue by contradiction. In particular, we assume  that $\widetilde{T}_0>0 $, otherwise there is nothing to prove. By  Proposition \ref{dual_ncob en},  system~(\ref{boun.sg2x2}) is null-controllable at time $T$ if and only if there exists $C>0$ such that any solution $\theta $ of the adjoint problem~(\ref{Padjoint}) satisfies the observability inequality: 
\begin{equation}
\Vert \theta (0)\Vert _{H^1_0(0,\pi ;\mathbb{R}^{2})}^{2}\leq
C\int_{0}^{T}\!\!|B^*D\theta _{x}(0,t)|^{2}\,dt,\quad \forall \theta
^{0}\in H^1_0\left( 0,\pi ;\mathbb{R}^{2}\right) .  \label{IOnug1}
\end{equation}
Let us recall that 
\begin{equation*}
\label{maxt1t2}
\widetilde{T}_0=max\left\lbrace T_1, \widetilde{T}_2\right\rbrace
\end{equation*}
with

\begin{equation*}
\label{t1t2}
T_1=\limsup_{k \longrightarrow +\infty}\frac{\log\Vert \psi_{1,k}\Vert_{H^1_0}}{k^2}, \quad 
\widetilde{T}_2=\limsup_{k \to +\infty } \frac{\log\frac{\left \Vert \psi_{2,k} -\psi_{1,i_k} \right\Vert}{ \left \vert \nu k^2- i_k^2 \right \vert}}{\nu k^2} =\limsup_{k \to +\infty } \frac{\log\frac{\left \Vert \psi_{2,k} -\psi_{1,i_k} \right\Vert}{ \left \vert \nu k^2- i_k^2 \right \vert}}{i_k^2}.
\end{equation*}

$\bullet $ Let us suppose   $\widetilde{T}_0=T_1$, and work with the particular solutions of the adjoint problem (\ref{Padjoint})  associated with initial data
$$\theta^0_{k} = \psi_{1,k} 
,$$ 
where $\psi_{1,k}$ is given by (\ref{psik}).
With this choice, the solution $\theta_{k}$ of the adjoint problem is given by 

$$
 \theta_{k} (x,t)=\psi_{1,k}(x) e^{-k^2(T-t)}.
$$
The observability inequality reads as 
 
$$
\left \Vert \theta_{k} (\cdot,0)\right \Vert^2_{H^1_0} \leq
C\int_{0}^{T}\left \vert B^{\ast }D \frac{\partial \theta_{k} (0,t)}{\partial x} (0,t)\right \vert^{2}\,dt.
$$
From the definition of $\widetilde{T}_1$, there exists a sub-sequence  $\{k_{n}\}_{n\ge1}$
 such that : 
$$
\widetilde{T}_1=\lim_{n\to+\infty}\frac{\log\Vert \psi_{1,k_n}\Vert_{H^1_0}}{k_n^2}\in]0,+\infty].
$$ 
in this case, for every $\varepsilon>0$, there exists a positive integer $n_{\varepsilon}\geq1$ such that 
$$ 
\Vert \psi_{1,k_n}\Vert_{H^1_0}>e^{k_n^2(T_1-\varepsilon)},\quad\forall n\ge n_{\varepsilon},
 $$
thus

$$
 \begin{cases}
 \begin{aligned}
&\Vert \theta_{k_n} (\cdot,0)\Vert_{H^1_0}^2=
e^{-2k_n^2T}\left \Vert \psi_{1,k_n} \right\Vert_{H^1_0}^2 \ge e^{-2k_n^2\left[ T-T_1 +\varepsilon\right] }  \underset{n\rightarrow +\infty}{\longrightarrow} +\infty,\\
&\int_{0}^{T}|B^{\ast }D \frac{\partial \theta_{k_n} (0,t)}{\partial x} (0,t)|^{2}\,dt=\int_{0}^{T} e^{-2k_n^2t}\,dt \underset{n\rightarrow +\infty}{\longrightarrow} 0.
\end{aligned}
\end{cases}
$$
This proves that the observability inequality  does not hold

$\bullet$  Let us suppose now, $\widetilde{T}_0=\widetilde{T}_2$ . Let us fix $ k \ge 1$ and work with the particular solutions associated with initial  data 
$$\theta^0_{k} = \psi_{1,i_k} - \psi_{2,k}
,$$ 
%where $\psi_{1,k}$ and $\psi _{2,j_k}$ are given by (\ref{psiik}).
with this choice, the solution $\theta_{k}$ of adjoint problem is given by 

$$
 \theta_{k} (x,t)=\psi_{1,i_k}(x) e^{-i_k^2(T-t)}- \psi_{2,k}(x)e^{-\nu k^2(T-t)}.
$$
We have
\begin{equation*}
 \begin{cases}
 \begin{aligned}
&\left \Vert \theta_{k} (\cdot,0)\right \Vert_{H^1_0} ^2=
\left \Vert e^{-i_k^2T} \left( \psi_{1,i_k} -\psi_{2,k} \right)
+
\psi_{2,k}  \left (e^{-i_k^2T}-  e^{-\nu k^2T}\right ) \right\Vert_{H^1_0}^2,
\\
&\int_{0}^{T}\left \vert B^{\ast }D \frac{\partial \theta_{k} (0,t)}{\partial x} (0,t)\right \vert^{2}\,dt=\int_{0}^{T}\left\vert   e^{-i_k^2t}-  e^{-\nu k^2t}\right\vert ^{2}\,dt,
\end{aligned}
\end{cases}
 %\label{NI=}
\end{equation*}
thus we obtain
\begin{equation*}
 \begin{cases}
&\Vert \theta_{k} (\cdot,0)\Vert_{H^1_0} \ge
e^{-2i_k^2T} \left \Vert \psi_{1,i_k} -\psi_{2,k} \right\Vert^2 -2e^{-i_k^2T} \left \Vert \psi_{1,i_k} -\psi_{2,k} \right\Vert \left \Vert \psi_{2,k} \right\Vert \left \vert i_k^2-\nu k^2\right \vert\\
\\
&\int_{0}^{T}|B^{\ast }D \frac{\partial \theta_{k} (0,t)}{\partial x} (0,t)|^{2}\,dt \leq   \left \vert i_k^2-\nu k^2 \right \vert^2,
\end{cases}
% \label{NI+}
\end{equation*}
and %from (\ref{NI=}) :
\begin{equation*}
 %\begin{array}{ccc}
\frac{\Vert \theta_{k}  (\cdot,0)\Vert^2_{H^1_0}  }{\int_{0}^{T}|B^{\ast }D \frac{\partial \theta_{k}  (0,t)}{\partial x} (0,t)|^{2}\,dt}  
\ge e^{-2i_k^2T} \frac{\left \Vert \psi_{1,i_k} -\psi_{2,k} \right\Vert^2}{ \left \vert i_k^2-\nu k^2 \right \vert^2} \left(1- 2e^{i_k^2T}\frac{\left \vert i_k^2-\nu k^2 \right \vert }{\left \Vert \psi_{1,i_k} -\psi_{2,k} \right\Vert} \left \Vert \psi_{2,i_k} \right\Vert \right).
%\end{array}
 %\label{Rkj}
\end{equation*}
From the definition of $\widetilde{T}_2$, there exists an increasing unbounded sub-sequence
 $\left(k_n\right)_{n \ge 1}$ such that
$$
\widetilde{T}_2=\lim_{n\to+\infty}\frac{\log\frac{\left \Vert \psi_{2,k_n} -\psi_{1,i_{k_n}} \right\Vert}{|\nu k_n^2- i_{k_n}^2|}}{i_{k_n^2}}\in]0,+\infty],
$$ 
in this case, for every $\varepsilon>0$, there exists a positive integer $n_{\varepsilon}\geq1$ such that 
$$ 
\frac{\left \Vert \psi_{2,k_n} -\psi_{1,i_{k_n}} \right\Vert}{|\nu k_n^2- i_{k_n}^2|}>e^{i_{k_n}^2(\widetilde{T}_2-\varepsilon)},\quad\forall n\ge n_{\varepsilon},
 $$
thus
 
 $$
 e^{-i_{k_n}^2T}\frac{\left \Vert \psi_{2,k_n} -\psi_{1,i_{k_n}} \right\Vert}{|\nu k_n^2- i_{k_n}^2|}>e^{i_{k_n}^2(\widetilde{T}_2-T-\varepsilon)} \underset{n\rightarrow +\infty}{\longrightarrow} +\infty, \quad \forall \varepsilon\in (0, \widetilde{T}_2-T).
 $$
 This proves that the observability inequality  does not hold and finishes the proof of the negative null controllability result of (\ref{boun.sg2x2}).

%%%%%%%%%%%%%%%% sous section nu>1 irrationnel %%%%%%%%%%%%%%
\subsection{Proof of Theorem \ref{Thm} : second case (b)}
%%%%%%%%%%%%%%%%%%%%%%%%%%%%%%%%%%%%%%%%%%%%%%%%%%%%%%%%%
In this subsection, our objective is to prove that system  (\ref{boun.sg2x2}) is null controllable  at time $T$  if $T > \widetilde{T}_0\in [0,\infty)$, when $\sqrt{\nu}<1$ . Thus, in the sequel we consider $\sqrt{\nu}<1$.

%%%%%%%%%%%%%%%%%%%%%%%%%%%%%%%%%%%%%%%%%%%%%%%%%
\subsubsection  {Positive null controllability result}
%%%%%%%%%%%%%%%%%%
 Let us start with the following crucial result.
\begin{prpstn}
\label{eppnew}
Let  $\sqrt{\nu}<1$. Let us define  the function
$$
\begin{array}{ccc}
J : &\mathbb{N}^{\ast}&\longrightarrow \mathbb{N}^{\ast} \\
&k& \longmapsto j_k,
\end{array}
$$
where, for any $k \in \mathbb{N}^*$, $j_k$ is the nearest integer to  $\frac{k}{\sqrt{\nu}}$  i.e 
$
\left \vert \frac{k}{\sqrt{\nu}} -j_k   \right \vert<\frac{1}{2} .
$
Thus for any $k \in \mathbb{N}^*$,
$$ \left \vert \frac{k}{\sqrt{\nu}}-j \right\vert >\frac{1}{2},\quad  \forall \,j \in \mathbb{N}^*,\, j\neq j_k .$$
Then
\begin{enumerate}

\item
The function $J$ is injective.

\item
$
 \widehat{J}=\mathbb{N}^*\setminus J(\mathbb{N}^*)=\{ \widehat{j}_k: k\ge 1\} 
$
is a infinite set,  where the elements of $\widehat{J}$ are classified in ascending order.
\end{enumerate}
\end{prpstn}
\noindent
\textit {{Proof of Proposition \ref{eppnew}}}.
The proof can be done in the same way as in Proposition \ref{epp}.
$\hfill \square $

 Thus, thanks to Proposition \ref{eppnew} we can reformulate (\ref{momInternnuir}). We say that the null controllability property at time $T$ for system (\ref{boun.sg2x2}) is equivalent to find $u\in L^2(0, T)$ such that:
 
\begin{equation}
 \label{momInternnuirbisnup1}
	\left\{ 	
	\begin{aligned}
&	 \int_0^T e^{-k^2t} u(T-t) \, dt=-e^{-k^2T}\langle y^0,\psi_{1,k} \rangle_{-1,i},\\
	%\left\{ 	
	%\begin{array}{l}
&\int_0^T e^{-\nu j_k^2t} u(T-t)\, dt=-e^{-\nu j_k^2T}\langle y^0,\psi_{2,j_k} \rangle_{-1,1}
 ,  \\ 
&\int_0^T e^{-\nu \widehat{j}_k^2t} u(T-t) \, dt=-e^{-\nu \widehat{j}_k^2T}\langle y^0,\psi_{2,\widehat{j}_k} \rangle_{-1,1},\quad  k\ge 1.
%\end{array}
%\right. 
\end{aligned}
\right.  
\end{equation}

We can now state our following main result.

\begin{prpstn}\label{control+}
There exists $\{q_k \}_{k \ge 1} \subset E_T$  such that 
\begin{enumerate}
\item
\begin{equation*}
 %\label{biorth}
 \left\{ 
 \begin{aligned}
 & \int_0^T e^{-k^2t} q_{k}(t) \, dt=-e^{-k^2T}\langle y^0,\psi_{1,k} \rangle_{-1,i},\\
&\int_0^T e^{-\nu j_k^2t} q_{k}(t)\, dt=-e^{-\nu j_k^2T}\langle y^0,\psi_{1,  j_k} \rangle_{-1,1}
 ,  \\ 
&\int_0^T e^{-\nu \widehat{j}_k^2t} q_{k}(t) \, dt=-e^{-\nu \widehat{j}_k^2T}\langle y^0,\psi_{2,\widehat{j}_k} \rangle_{-1,1},\quad  k\ge 1,
\\
&\int_0^T  e^{-k^2t}q_{i}(t)dt=\int_0^T  e^{-\nu j_k^2t}q_{i}(t)dt=\int_0^T  e^{-\nu \widehat{j}_k^2 t}q_{i}(t)dt=0, \, k\neq i.
    \end{aligned}
 \right.
\end{equation*}

\item If $T>\widehat{T}_0=\max\{T_1, \widehat{T}_2\}$ (see (\ref{T0hat}) then we infer that an explicit solution $u$  of moment problem (\ref{momInternnuirbisnup1})  given by

\begin{equation*}
% \label{unuir}
u(t)=\sum_{k\ge 1}q_k(T-t).
\end{equation*}
\end{enumerate}

\end{prpstn}

\bigskip

\noindent
\textit {{Proof of Proposition \ref{control+}}}.
Just use here the same reasoning as in the previous Proposition \ref{control}. 
$\hfill \square $
%%%%%%%%%%%%%%%%
\subsubsection  {The negative null controllability result }
%%%%%%%%%%%%%%
Let us prove that if $0<T<\widehat{T}_0$, then system (\ref{boun.sg2x2}) is not null controllable at time $T$. We argue by contradiction. Let us first recall that 

\begin{equation*}
%\label{maxt1t2 nu2}
\widehat{T}_0=max\left\lbrace T_1, \widehat{T}_2\right\rbrace
\end{equation*}
where

\begin{equation*}
%\label{t1t2 nu2}
T_1=\limsup_{k \longrightarrow +\infty}\frac{\log\Vert \psi_{1,k}\Vert_{H^1_0}}{k^2} \quad \hbox{and}\quad
\widehat{T}_2=\limsup_{k \to +\infty } \frac{\log\frac{\left \Vert \psi_{1,k} -\psi_{2,j_k} \right\Vert}{ \left \vert k^2- \nu j_k^2 \right \vert}}{k^2}   =\limsup_{k \to +\infty } \frac{\log\frac{\left \Vert \psi_{1,k} -\psi_{1,j_k} \right\Vert}{ \left \vert k^2- \nu j_k^2 \right \vert}}{\nu j_k^2}.
\end{equation*}
Assume  that system ~\eqref{boun.sg2x2}  is  null-controllable at time $T<\widehat{T}_0 $. System~(\ref{boun.sg2x2}) is null-controllable at time $T$ if and only if there exists $C>0$ such that any solution $\theta $ of the adjoint problem~(\ref{Padjoint}) satisfies the observability inequality: 
\begin{equation*}
\Vert \theta (0)\Vert _{H^1_0(0,\pi ;\mathbb{R}^{2})}^{2}\leq
C\int_{0}^{T}\!\!|B^*D\theta _{x}(0,t)|^{2}\,dt,\quad \forall \theta
^{0}\in H^1_0\left( 0,\pi ;\mathbb{R}^{2}\right) . 
%\label{IOnug1 nu2}
\end{equation*}
$\bullet$ Let us suppose   $\widehat{T}_0=T_1$  and work with the particular solutions of adjoint problem (\ref{Padjoint})  associated with initial data
$$\theta^0_{k} = \psi_{1,k} 
,$$ 
where $\psi_{1,k}$ is given by (\ref{psik}).
With this choice, the solution $\theta_{k}$ of adjoint problem is given by 

$$
 \theta_{k} (x,t)=\psi_{1,k}(x) e^{-k^2(T-t)}.
 $$
The observability inequality reads as 
 
$$
\left \Vert \theta_{k} (\cdot,0)\right \Vert^2_{H^1_0} \leq
C\int_{0}^{T}\left \vert B^{\ast }D \frac{\partial \theta_{k} (0,t)}{\partial x} (0,t)\right \vert^{2}\,dt.
$$

From the definition of $T_1$,  there exists a subsequence  $\{k_{n}\}_{n\ge1}$
 such that : 
$$
T_1=\lim_{n\to+\infty}\frac{\log\Vert \psi_{1,k_n}\Vert_{H^1_0}}{k_n^2}\in]0,+\infty],
$$ 
in this case, for every $\varepsilon>0$, there exits a positive integer $n_{\varepsilon}\geq1$ such that 
$$ 
\Vert \psi_{1,k_n}\Vert_{H^1_0}>e^{k_n^2(T_{1}-\varepsilon)},\quad\forall n\ge n_{\varepsilon},
 $$
thus, 
\begin{equation*}
 \begin{cases}
 \begin{aligned}
&\Vert \theta_{k_n} (\cdot,0)\Vert_{H^1_0}^2=
e^{-2k_n^2T}\left \Vert \psi_{1,k_n} \right\Vert_{H^1_0}^2 \ge e^{-2k_n^2\left[ T-T_1 +\varepsilon\right] }\longrightarrow +\infty,\quad \hbox{when}\;n\longrightarrow +\infty,
  \\
&\int_{0}^{T}|B^{\ast }D \frac{\partial \theta_{k_n} (0,t)}{\partial x} (0,t)|^{2}\,dt=\int_{0}^{T} e^{-2k_n^2t}\,dt \longrightarrow 0,\quad \hbox{when}\; n\longrightarrow +\infty.
 \end{aligned}
\end{cases}
% \label{NI1}
\end{equation*}
This proves that the observability inequality  does not hold.

$\bullet$ Assume $\widehat{T}_0=\widehat{T}_2$ . Let us fix $ k \ge 1$ and work with the particular solutions associated with initial  data 
$$\theta^0_{k} = \psi_{1,k} - \psi_{2,j_k}
.$$ 
%where $\psi_{1,k}$ and $\psi _{2,j_k}$ are given by (\ref{psiik}).
With this choice, the solution $\theta_{k}$ of adjoint problem is given by 
$$
 \theta_{k} (x,t)=\psi_{1,k}(x) e^{-k^2(T-t)}- \psi_{2,j_k}(x)e^{-\nu j_k^2(T-t)}.
$$
We have
\begin{equation*}
 \begin{cases}
 \begin{aligned}
&\left \Vert \theta_{k} (\cdot,0)\right \Vert_{H^1_0} ^2=
\left \Vert e^{-k^2T} \left( \psi_{1,k} -\psi_{2,j_k} \right)
+
\psi_{2,j_k}  \left (e^{-k^2T}-  e^{-\nu j_k^2T}\right ) \right\Vert_{H^1_0}^2, 
  \\
&\int_{0}^{T}\left \vert B^{\ast }D \frac{\partial \theta_{k} (0,t)}{\partial x} (0,t)\right \vert^{2}\,dt=\int_{0}^{T}\left\vert   e^{-k^2t}-  e^{-\nu j_k^2t}\right\vert ^{2}\,dt.
\end{aligned}
\end{cases}
 \label{NI=}
\end{equation*}
Thus we obtain
\begin{equation*}
 \begin{cases}
 \begin{aligned}
&\Vert \theta_{k} (\cdot,0)\Vert_{H^1_0} \ge
e^{-2k^2T} \left \Vert \psi_{1,k} -\psi_{2,j_k} \right\Vert^2 -2e^{-k^2T} \left \Vert \psi_{1,k} -\psi_{2,j_k} \right\Vert \left \Vert \psi_{2,j_k} \right\Vert \left \vert k^2-\nu j_k^2\right \vert,\\
&\int_{0}^{T}\left \vert B^{\ast }D \frac{\partial \theta_{k} (0,t)}{\partial x} (0,t)\right \vert^{2}\,dt \leq   \left \vert k^2-\nu j_k^2 \right \vert^2,
\end{aligned}
\end{cases}
 %\label{NI+}
\end{equation*}
and %from (\ref{NI=}) :
\begin{equation*}
 %\begin{array}{ccc}
\frac{\Vert \theta_{k}  (\cdot,0)\Vert^2_{H^1_0}  }{\displaystyle\int_{0}^{T}\left \vert B^{\ast }D \frac{\partial \theta_{k}  (0,t)}{\partial x} (0,t)\right\vert^{2}\,dt}  
\ge e^{-2k^2T} \frac{\left \Vert \psi_{1,k} -\psi_{2,j_k} \right\Vert^2}{ \left \vert k^2-\nu j_k^2 \right \vert^2} \left(1- 2e^{k^2T}\frac{\left \vert k^2-\nu j_k^2 \right \vert }{\left \Vert \psi_{1,k} -\psi_{2,j_k} \right\Vert} \left \Vert \psi_{2,j_k} \right\Vert \right).
%\end{array}
% \label{Rkj}
\end{equation*}
From the definition of $\widehat{T}_2$, there exists an increasing unbounded subsequence
 $\left(k_n\right)_{n \ge 1}$ such that
$$
\widehat{T}_2=\lim_{n\to+\infty}\frac{\log\frac{\left \Vert \psi_{1,k_n} -\psi_{2,j_{k_n}} \right\Vert}{|k_n^2- \nu j_{k_n}^2|}}{k_n^2}\in]0,+\infty],
$$ 
in this case, for every $\varepsilon>0$, there exits a positive integer $n_{\varepsilon}\geq1$ such that 
$$ 
\frac{\left \Vert \psi_{1,k_n} -\psi_{2,j_{k_n}} \right\Vert}{|k_n^2- \nu j_{k_n}^2|}>e^{k_n^2(\widehat{T}_2-\varepsilon)},\quad\forall n\ge n_{\varepsilon},
 $$
thus, 
 
 $$
 e^{-k_n^2T}\frac{\left \Vert \psi_{2,k_n} -\psi_{1,i_{k_n}} \right\Vert}{|\nu k_n^2- i_{k_n}^2|}>e^{k_n^2(\widehat{T}_2-T-\varepsilon)} \longrightarrow +\infty, \quad \forall \varepsilon\in (0, \widehat{T}_2-T).
 $$
 This proves that the observability inequality  does not hold and finishes the proof of the negative null controllability result of (\ref{boun.sg2x2}).
 
 %%%%%%%%%%%%%%%%%%%%%%%%%%%%%%%%%%%%%%%%%%%%%%%%%%%%%%%%%%%%%%%%%%%%%%%%%
\section{ Null-controllability result with an rational diffusion coefficient}
%%%%%%%%%%%%%%%%%%%%%%%%%%%%%%%%%%%%%%%%%%%%%%%%%%%%%%%%%%%%%%

%%%%%%%%%%%%%%%%%%%%%%%%%%%%%%%%%%%%%%%%%%%%%%%%%%%%%%%%%%%
\subsection{ Some preliminary results}
\label{spr1nur}
%%%%%%%%%%%%%%%%%%%%%%%%%%%%%%%%%%%%%%%%%%%%%%%%%%%%%%%%%%%%%%%

In this subsection we will give some properties which will be used below. We are interested in studying the spectrum of the operators $L^{\ast}$. Let us define the sets

  \begin{equation} 
\left\{ 
	\begin{aligned}
	& \Lambda_1:=\left\{ \frac{i_0^2}{j_0^2}j^2:j\in \mathbb{N}^{\ast}\setminus j_0\mathbb{N}^{\ast}\right\} ,\\
	 &\Lambda_2:=\left\{ k^2: k\in \mathbb{N}^{\ast} \setminus i_0\mathbb{N}^{\ast}  \right\},\\
	& \Lambda_3:=\left\{ i_0^2l^2:l\ge 1\right\}.
\end{aligned}
\right.
\label{Lambdas}
\end{equation}
Observe that $\Lambda_1$, $\Lambda_2$ and $\Lambda_3$ are disjoint sets.
\begin{rmrk}
If $i_0=1$ then $\Lambda_2=\emptyset$, $\Lambda_1:=\{ \frac{j^2}{j_0^2}:j\in \mathbb{N}^{\ast}\setminus j_0\mathbb{N}^{\ast}\}$ and $\Lambda_3:=\{ l^2:l\ge 1\}$. Thus in the sequel, the controllability of the system (\ref{boun.sg2x2}) will be studied in only when  $\sqrt{\nu}=\frac{i_0}{j_0} \in \mathbb{Q}_{+}^*$ and $i_0>1$. The case $i_0=1$ is much easier and it is left to the reader.
\end{rmrk}
In the sequel, for the proof of Theorem \ref{Thm+}, we shall use the following notations
\begin{equation*}
 I(\zeta)=\int_{0}^{\pi}q(s)\varphi_{\sqrt{\zeta}}(s)\sqrt{\frac{\pi}{2}}\sin\left(\sqrt{\frac{\zeta}{\nu}}(\pi-s)\right)ds, \quad \zeta \in \Lambda_2 \cup \Lambda_3=\{ k^2:k\ge 1\} ,
% \label{Ik}
\end{equation*}
%where $\sigma(L^{\ast }) $ is the spectrum of $L^{\ast }$, and
and $T_0(\nu,q)$ becomes
\begin{equation} 
\label{T0alpha}
T_0(\nu,q):=\limsup_{\zeta \in \Lambda_2 \cup \Lambda_3,\;\zeta \to +\infty }  \frac{-\log  \left \vert I(\zeta)  \right \vert }{ \zeta} \in [0, +\infty].
\end{equation}

\begin{prpstn}\label{fam_nu_rnur1}
Let $A_0$ be given by ~(\ref{matDAB}) and consider the operator $(L,D(L))$ given by ~(\ref{opL}) and its
adjoint $L^{\ast }$. 

\begin{enumerate}

\item 
The spectrum of $L^{\ast }$ is given by 
\begin{equation}
\sigma(L^{\ast })=\Lambda_1 \cup \Lambda_2 \cup \Lambda_3.
\label{sigmaLast}
\end{equation}

\item
 Given $j\in \mathbb{N}^{\ast}\setminus j_0\mathbb{N}^{\ast}$,  $k\in \mathbb{N}^{\ast}\setminus i_0\mathbb{N}^{\ast}$  and  $l\ge 1$, let us introduce
 
 \begin{equation} 
 \Phi _{2,j}^{\ast }:=\left( 
\begin{aligned}
0 \\ 
\varphi _{j}
\end{aligned}
\right),
\,
\widetilde{\Phi}_{1,k}^{\ast }:=\left( 
\begin{aligned}
\varphi _{k} \\ 
\widetilde{\psi} _{k}%
\end{aligned}%
\right)\,\hbox{ and}\quad
\widehat{\Phi} _{1,l}^{\ast }:=\left( 
\begin{aligned}
\varphi _{i_0l} \\ 
\widehat{\psi} _{l}%
\end{aligned}%
\right).
\label{vectprop}
\end{equation} 
\begin{enumerate}

\item Given $j\in \mathbb{N}^{\ast}\setminus j_0\mathbb{N}^{\ast}$, then 

\begin{equation*}
\left( L^{\ast }-\nu j^2 I_{d}\right) \Phi _{2,j}^{\ast }=0.
%\label{phi2j}
\end{equation*}

\item Given $k\in \mathbb{N}^{\ast}\setminus i_0\mathbb{N}^{\ast}$, then 

\begin{equation*}
 \left(L^{\ast }-k^2 I_{d}\right) \widetilde{\Phi}_{1,k}^{\ast }=0.
%\label{phitilde1k}
\end{equation*}

where, $\widetilde{\psi}_k$  is the unique solution of following problem:

\begin{equation} 
 \label{psikeq}
\left\{ 
	\begin{aligned}
& - \widetilde{\psi}_k-\frac{k^{2}}{\nu} \widetilde{\psi}_k=-\frac{1}{\nu}q\varphi_{k} \quad \hbox{in }(0,\pi ), \\ 
&\widetilde{\psi}_k(0)= \widetilde{\psi}_k(\pi )=0. \\
\end{aligned}%
\right.
\end{equation}%
  Moreover, an explicit expression of $\widetilde{\psi}_{k}$
  is given by :
\begin{equation} 
 \label{tpsiksol}
\widetilde{\psi}_{k}(x)=\widetilde{\psi}'_{k}(0)\frac{\sqrt{\nu}}{k}\sin\left(\frac{kx}{\sqrt{\nu}}\right)+\frac{\sqrt{\nu}}{\nu k}\int_{0}^{x}\sin\left(\frac{k}{\sqrt{\nu}}(x-\xi)\right)q(\xi)\varphi_{k}(\xi)\,d\xi,
\end{equation}%
with (see Remark \ref{remark} )

\begin{equation} 
 \label{psi'tildek0}
\widetilde{\psi}'_{k}(0)=-\frac{I(k^2)}{\nu \sqrt{\frac{\pi}{2}}\sin\left(\frac{k\pi}{\sqrt{\nu}}\right)}.
\end{equation}%
\item Given $l\ge 1$, then 

\begin{equation*}
 \left\{ 
	\begin{aligned}
&	\left(L^{\ast }-i_0^2l^2 I_{d}\right) \Phi _{2,j_0l}^{\ast }=0,\\
&	 \left(L^{\ast }-i_0^2l^2 I_{d}\right) \widehat{\Phi} _{1,l}^{\ast }=I (i_0^2l^2)\Phi _{2,j_0l}^{\ast },\quad \forall l \ge 1,
	 \end{aligned}
	 \right.
%\label{phihat1k}
\end{equation*}

and $\widehat{\psi}_l$  is the unique solution of following problem:

\begin{equation*} 
 %\label{psikeq}
\left\{ 
	\begin{aligned}
 &- \widehat{\psi}_l-j_0^2l^{2} \widehat{\psi}_l=\frac{1}{\nu}[I(i_0^2 l^2)\varphi_{j_0l}-q\varphi_{i_0l}] \quad \hbox{in  }(0,\pi ), \\ 
&\widehat{\psi}_l(0)= \widehat{\psi}_l(\pi )=0, \\
&\int_{0}^{\pi}\widehat{\psi}_l(x) \varphi_{j_0l}(x)\,dx=0,
\end{aligned}%
\right.
\end{equation*}%
 an explicit expression of $\widehat{\psi}_{l}$
  is given by :
\begin{equation} 
 \label{hpsiksol-}
  \widehat{\psi}_{l}(x)=\alpha_l \varphi_{j_0l}(x)+\beta_l(x), \quad \hbox{for all}\quad x\in (0,\pi),
 \end{equation}%
  where
  
\begin{equation*} 
% \label{hpsiksol}
 \left\{ 
	\begin{aligned}
\beta_l(x)&=-\frac{1}{\nu j_0l}\int_{0}^{x}\sin\left(j_0l(x-\xi)\right)[I (i_0^2l^2)\varphi_{j_0l}(\xi)-q(\xi)\varphi_{i_0l}(\xi)]\,d\xi \\
\alpha_l& =\frac{1}{\nu j_0l}\int_{0}^{\pi}\int_{0}^{x}\sin\left(j_0l(x-\xi)\right)[I (i_0^2l^2)\varphi_{j_0l}(\xi)-q(\xi)\varphi_{i_0l}(\xi)]\varphi_{j_0l}(x)\,d\xi\,dx.
\end{aligned}%
\right.
\end{equation*}%
Moreover (see Remark \ref{remark} for instance)

\begin{equation} 
 \label{beta'l0}
\beta'_l(x)=-\frac{1}{\nu }\int_{0}^{x}\cos\left(j_0l(x-\xi)\right)[I (i_0^2l^2)\varphi_{j_0l}(\xi)-q(\xi)\varphi_{i_0l}(\xi)]\,d\xi,
\end{equation}%
thus  $\beta'_l(0)=0$.
\end{enumerate}
  \end{enumerate}
  \end{prpstn}
  
  \bigskip

In the next result we are going to give some properties of  $\widetilde{\psi}_k$ ( see \ref{tpsiksol})   and $\widehat{\psi}_l$ (see \ref{hpsiksol-}). These properties will be used later and will be crucial in the proof of Theorem \ref{Thm+}.

\smallskip

\begin{prpstn}\label{prop2est}
\begin{enumerate}
\item Let us fix $q\in L^{\infty}(0,\pi)$ and take $k \in \mathbb{N}^{\ast}\setminus i_0\mathbb{N}^{\ast}$. Then, one has:
\begin{equation} 
 \label{tpsiksolprime}
\widetilde{\psi}'_{k}(x)=-\frac{I(k^2)}{\nu \sqrt{\frac{\pi}{2}}\sin\left(\frac{k\pi}{\sqrt{\nu}}\right)}\cos\left(\frac{kx}{\sqrt{\nu}}\right)+\frac{1}{\nu }\int_{0}^{x}\cos\left(\frac{k}{\sqrt{\nu}}(x-\xi)\right)q(\xi)\varphi_{k}(\xi)\,d\xi.
\end{equation}%
In addition, there exists a constant $C > 0$ such that
\begin{equation*} 
\left \Vert \widetilde{\psi}'_k\right \Vert_{L^{\infty}(0,\pi)} \le C,\quad k \in \mathbb{N}^{\ast}\setminus i_0\mathbb{N}^{\ast}.
%\label{tpsi'kest}
\end{equation*}
\item Let us take $l \ge 1$. There exists a constant $C > 0$ such that
\begin{equation} 
\left \vert \alpha_l \right\vert \le C \quad \hbox{and}\quad \left \Vert \widehat{\psi}'_l \right \Vert_{L^{\infty}(0,\pi)} \le C, \quad l\ge 1.
\label{hpsi'kest}
\end{equation}
\end{enumerate}
\end{prpstn}

\noindent
\textit {Proof of Proposition \ref{prop2est}}. 

\begin{enumerate}
\item Let us fix $k \in \mathbb{N}^{\ast}\setminus i_0\mathbb{N}^{\ast}$, the expression (\ref{tpsiksolprime}) can be deduced from (\ref{y'gene}) and (\ref{cnsygene}). Moreover,  there exists a constant $C > 0$ such that %the property (\ref{tpsi'kest}) can be easily deduced from (\ref{tpsiksolprime}). 

$$
\begin{aligned}
\left \Vert \widetilde{\psi}'_k\right \Vert_{L^{\infty}(0,\pi)} &\le  \frac{C}{\left \vert \sin\left( \frac{k \pi}{\sqrt{\nu}}\right) \right \vert }+C, \\
&= \frac{C}{\left \vert \sin\left( \pi (n_k - \frac{k j_0 }{i_0}) \right) \right \vert }+ C  & &\left(\hbox{for each $k$,}\; n_k \; \hbox{is the nearest integer of} \; \frac{kj_0}{i_0}\right),\\
 &\le  \frac{\pi C}{2\left \vert  \pi (n_k - \frac{k j_0 }{i_0}) \right \vert }+C  & &\left( 0<\left \vert   n_k - \frac{k j_0 }{i_0} \right \vert \le\frac{1}{2},\: \forall \, x\in [0,\frac{\pi}{2}],\: \frac{2}{\pi}x \le \sin x\right),\\
 & \le  \frac{i_0 C}{2 }+C  & &\left(\hbox{ since } 0<\frac{1}{i_0} \le\left \vert   n_k - \frac{k j_0 }{i_0} \right \vert \le \frac{1}{2} \right).\\
\end{aligned}
$$
\item The properties (\ref{hpsi'kest}) can be deduced from (\ref{hpsiksol-}) . This finalizes the proof.
\end{enumerate}
$\hfill \square $

\smallskip
  
\noindent
\textit {Proof of Proposition \ref{fam_nu_rnur1}}. 
%Let us assume $\sqrt{\nu}=\frac{i_0}{j_0}  \in \mathbb{Q}_{+}^*$ and $k,j,l\ge 1$.

First, let $ \lambda $ be an eigenvalue of $L^{\ast}$ and $y=(y_1,y_2)^{T}$ an associated eigenfunction. Thus $y$ is a solution of the following problem:

\begin{equation*} 
\left\{ 
	\begin{aligned}
&	 -y''_1=\lambda y_1 \: &\hbox{ in }\:(0,\pi ) \\ 
 &qy_1-\nu y''_2=\lambda y_2 \: &\hbox{ in }\: (0,\pi) \\ 
&y_1(0)=y_2(0)=0,\quad y_1(\pi)=y_2(\pi)=0.\\
\end{aligned}%
\right.  
\end{equation*}%
\begin{itemize}

\item If $y_1\equiv 0$, then,  $\lambda=\nu j^2$ is an eigenvalue of $L^{\ast}$ and taking $y_2=\varphi_{j}$, we obtain $\Phi^{\ast}_{2,j}$ as associated eigenfunction of $L^{\ast}$. 

\item Now assume that $y_1\not\equiv 0$, then $\lambda=k^2$ and $y_1 = \varphi_{k}$ is a solution to the first o.d.e. Inserting this expression in the second equation, we get for $y_2$ :
\begin{equation} 
\label{y2eq}
\left\{ 
	\begin{array}{ll}
y''_2+\frac{k^2}{\nu} y_2 =\frac{1}{\nu} q \varphi_{k} \: \hbox{in }\: (0,\pi) \\ 
	\noalign{\smallskip}	
y_2(0)=y_2(\pi)=0.\\
\end{array}%
\right.  
\end{equation}%
Let us multiply (\ref{y2eq}) by $\sin\left(\frac{k}{\sqrt{\nu}}(\pi-s)\right)$ and using partial integration, we obtain :
%A necessary and sufficient condition for the previous nonhomogeneous Sturm-Liouville problem to have a solution is that
\begin{equation*} 
%\label{cnsy2}
\left\{ 
	\begin{aligned}
& y_{2}'(0)=\frac{-\int_{0}^{\pi}q(s)\varphi_{k}(s)\sin\left(\frac{k}{\sqrt{\nu}}(\pi-s)\right)ds}{\nu \sin\left(\frac{k\pi}{\sqrt{\nu}}\right)}=-\frac{I (k^2) }{\nu \sqrt{\frac{\pi}{2}}\sin\left(\frac{k\pi}{\sqrt{\nu}}\right)}, &\hbox{if}&\: k \neq i_0l, \\
&  \int_{0}^{\pi}q(s)\varphi_{i_0 l}(s)\varphi_{j_0l}(s)ds=I (i_0^2l^2)=0, &\hbox{if}&\: k = i_0l. \
  \end{aligned}%
\right. 
 \end{equation*}%
%\end{itemize} 

\begin{itemize}

\item Observe that, suppose $k=i_0l$, thus (\ref{y2eq}) admits a solution if and only if $I (i_0^2l^2)=0$. In this case, $i_0^2l^2$  is a double eigenvalue of $L^{\ast}$ (take $j=j_0l$).
  From the above considerations, it is clear that if $I (i_0^2l^2)\neq 0$, then the eigenvalue  $i_0^2l^2$ of $L^{\ast}$ is simple and  $\Phi^{\ast}_{2,j_0l}$  is an associated eigenfunction. Observe that, taking $\widehat{\Phi}^{\ast}_{1,l}=(z_1,z_2)^{T}$ the equation $ \left(L^{\ast }- i_0^2l^2I_{d}\right) \widehat{\Phi} _{1,l}^{\ast }=c\Phi_{2,j_0l}$ writes:
\begin{equation*} 
\left\{ 
	\begin{aligned}
	& -z''_1-i_0^2l^2 z_1=0 \: &\hbox{in }&(0,\pi ), \\ 
& -\nu z''_2-i_0^2l^2z_2 =[c\varphi_{j_0l}-qz_1] \:& \hbox{in }&(0,\pi), \\ 
&z_1(0)=z_2(0)=0,\quad z_1(\pi)=z_2(\pi)=0.\\
\end{aligned}%
\right.  
\end{equation*}%
Thus,  choosing $z_1 = \varphi_{i_0l}$ and inserting this expression in the second equation, we get :
\begin{equation} 
\label{z2}
\left\{ 
	\begin{aligned}
&- z''_2-j_0^2l^2 z_2 =\frac{1}{\nu}[c\varphi_{j_0l}-q\varphi_{i_0l}] \quad \hbox{in }\: (0,\pi), \\
&z_2(0)=z_2(\pi)=0.
\end{aligned}%
\right.  
\end{equation}%
A necessary and sufficient condition for the previous no-homogeneous Sturm-Liouville problem to
have a solution is that
\begin{equation*} 
\int _0^{\pi}[c\varphi^2_{j_0l}(x)-q(x)\varphi_{i_0l}(x)\varphi_{j_0l}(x)] dx=0 ,\quad \hbox{i.e} \quad c=\int_{0}^{\pi}q\varphi_{i_0l}\varphi_{j_0l}ds.
%\label{z2cns1}
\end{equation*}
With this value of c, (\ref{z2}) has a continuum of solutions. This proves that, for $k=i_0l$, $\widehat{\Phi}_{1,l}$ is a generalized eigenfunction of $L^{\ast}$ associated to $\lambda=i_0^2l^2$. \item Otherwise, if  $k\neq i_0l$, then
\begin{equation*} 
%\label{y'2(0)}
y_{2}'(0)=-\frac{\int_{0}^{\pi}q(s)\varphi_{k}(s)\sin\left(\frac{k}{\sqrt{\nu}}(\pi-s)\right)ds}{\nu \sin\left(\frac{k\pi}{\sqrt{\nu}}\right)}=-\frac{I(
k^2)}{\nu \sqrt{\frac{\pi}{2}}\sin\left(\frac{k\pi}{\sqrt{\nu}}\right)},
 \end{equation*}%
and $\widetilde{\Phi}_{1,k}$ is a  eigenfunction of $L^{\ast}$ associated to $\lambda=k^2$ with $k\neq i_0l$.
\end{itemize}
\end{itemize}
 $\hfill \square $
 
\begin{lmm}
\label{lb}
The family
 $\mathcal{B^{\ast}}=\left\{ \Phi^{\ast}_{2,j}:j \in \mathbb{N}^{\ast} \setminus j_0\mathbb{N}^{\ast} ;  \widetilde{\Phi}^{\ast} _{1,k}:k\in \mathbb{N}^{\ast} \setminus i_0\mathbb{N}^{\ast}; \Phi^{\ast}_{2,j_0l},\widehat{\Phi}^{\ast}_{1,l} :l \ge 1 \right \}$ 
 is complete in $L^{2}(0,\pi ;\mathbb{R}^2)$. 
\end{lmm}
\noindent
\textit {Proof of Lemma \ref{lb}}.
Indeed, if $f=(f_1,f_2)$ is such that

$$
\left\{  \begin{aligned}
&\langle f,\Phi_{2,j}^{\ast}\rangle=0,   &j& \in \mathbb{N}^{\ast} \setminus j_0\mathbb{N}^{\ast}, \\
&\langle f,\widetilde{\Phi}_{1,k}^{\ast}\rangle=0,   &k& \in \mathbb{N}^{\ast} \setminus i_0\mathbb{N}^{\ast},\\
&\langle f,\Phi_{2,j_0l}^{\ast}\rangle=0,   &l&\ge 1,\\
&\langle f,\widehat{\Phi}_{1,l}^{\ast}\rangle=0,  &l&\ge 1.
\end{aligned} \right.
$$
Then in particular
$$
\left\{  \begin{aligned}
&\langle f_2,\varphi_{j}\rangle=0, &j& \in \mathbb{N}^{\ast} \setminus j_0\mathbb{N}^{\ast}, \\
&\langle f_1,\varphi_{k}\rangle+\langle f_2,\widetilde{\psi}_{k}\rangle =0, &k& \in \mathbb{N}^{\ast} \setminus i_0\mathbb{N}^{\ast},\\
&\langle f_2,\varphi_{j_0l}\rangle=0, &l&\ge 1,\\
&\langle f_1,\varphi_{i_0l}\rangle+\langle f_2,\widehat{\psi}_{l}\rangle =0, 
&l&\ge 1.
\end{aligned}
\right.
$$
This implies that $f_1=f_2=0$ (since $\{\varphi_{k} \}_{k\ge 1}$ is an orthonormal basis in $L^2(0,\pi)$) and proves the completeness of $\mathcal{B^{\ast}}$. 
$\hfill \square $
%%%%%%%%%%%%%%%%%%%%%%%%%%%%%%%%%%%%%%%%%%%%%%%%
%%%%%%%%%%%%%%%%%%%%%%%%%%%%%%%%%%%%%%%%%%%%%%%%%%%%%%%%%%%%%%%%%%%%%%%%%
\subsection{ Proof of Theorem \ref{Thm+}  }
\label{proofthm+}
%%%%%%%%%%%%%%%%%%%%%%%%%%%%%%%%%%%%%%%%%%%%%%%%%%%%%%%%%%%%%%
%%%%%%%%%%%%%%%%%%%%%%%%%%%%%%%%%
%In this section we will prove Theorem  (\ref{Thm+}). To this end, l

\subsubsection{ Approximate controllability}
This subsection is devoted to proving the approximate controllability of system (\ref{boun.sg2x2}), that is to say, the first point of Theorem \ref{Thm+}. To this end, we are going to apply the property (\ref{FAT66}).

\textbf{Necessary condition:} Let us suppose that condition (\ref{cndca+}) does not hold i.e. there exists $k_0 \in \mathbb{N}^{\ast} \setminus i_0\mathbb{N}^{\ast}$ (resp. $l_0 \in \mathbb{N}^{\ast} $) such that $I( k_0^2)=0$  (resp. $I( i_0^2l_0^2)=0$ ).
\begin{itemize}
\item If $I( k_0^2)=0$, with $k_0 \in \mathbb{N}^{\ast} \setminus i_0\mathbb{N}^{\ast}$, then we deduce that 
$$\widetilde{\Phi}_{1,k_0}^{\ast}=
\left( 
\begin{array}{l}
\varphi _{k_0}\\ 
\widetilde{\psi}_{k_0}
\end{array}
\right) \quad \hbox{(see (\ref{vectprop}))},
$$
is an non-trivial eigenfunction associated with the eigenvalue $k_0^2$ of the operator $L^{\ast}$ satisfying
$$
B^{\ast}D \frac{\partial \widetilde{\Phi}_{1,k_0}^{\ast}}{\partial x}(0)=\nu \widetilde{\psi}'_{k_0} (0)=-\frac{I(k_0^2)}{\sqrt{\frac{\pi}{2}}\sin\left(\frac{k_0\pi}{\sqrt{\nu}}\right)}=0\quad \hbox{(see (\ref{psi'tildek0}})).
$$

\item If  $I( i_0^2 l_0^2)=0$, with $l_0 \in \mathbb{N}^{\ast} $, then we deduce that 
$$\widehat{\Phi}_{1,l_0}^{\ast}-\alpha_{l_0} \Phi_{2,j_0 l_0}^{\ast}=\left( 
\begin{array}{l}
\varphi _{i_0 l_0}\\ 
\beta_{l_0}
\end{array}
\right)\quad \hbox{(see (\ref{vectprop}) and (\ref{hpsiksol-}))},
$$
is an non-trivial eigenfunction associated with the eigenvalue $i_0^2 l_0^2$ of the operator $L^{\ast}$ satisfying
$$
B^{\ast}D \frac{\partial \left( \widehat{\Phi}_{1,l_0}^{\ast}-\alpha_{l_0} \Phi_{2,j_0 l}^{\ast} \right)}{\partial x}(0)=\nu\beta'_{l_0} (0)=0\quad \hbox{(see (\ref{beta'l0}))}.
$$

\end{itemize}
Thus, system (\ref{boun.sg2x2}) is not approximately controllable at time $T$.

\textbf{Sufficient condition:} Let us suppose that condition (\ref{cndca+})  hold. The set of the eigenvectors associated with the eigenvalue $\nu k^2$  of $L^{\ast}$ is generated by $\Phi_{2,k}^{\ast}$
 (see Proposition \ref{fam_nu_ir}). In this case, we remark that for all $k\in \mathbb{N}^{\ast}$
 
 $$
B^{\ast}D \frac{\partial \Phi_{2,k}^{\ast}}{\partial x}(0)=\nu \varphi'_k(0)=k\nu \sqrt{\frac{2}{\pi}}\neq 0.
$$
Moreover, the set of the eigenvectors associated with the eigenvalue $k^2$, $k \in \mathbb{N}^{\ast} \setminus i_0\mathbb{N}^{\ast}$,  of $L^{\ast}$ is generated by $\widetilde{\Phi}_{1,k}^{\ast}$ (see Proposition \ref{fam_nu_ir}). In this case, we remark that for all $k \in \mathbb{N}^{\ast} \setminus i_0\mathbb{N}^{\ast}$
 
 $$
B^{\ast}D \frac{\partial \widetilde{\Phi}_{1,k}^{\ast}}{\partial x}(0)=\nu \widetilde{\psi}'_k(0)=-\frac{I(k^2)}{\sqrt{\frac{\pi}{2}}\sin\left(\frac{k\pi}{\sqrt{\nu}}\right)}\neq 0\quad \hbox{(see (\ref{psi'tildek0}))}.
$$
We conclude with the help of (\ref{FAT66}).

%$\hfill \square$

%%%%%%%%%%%%%%%%%%%%%%%%%%%%%%%%%%%%%%%%%%%%%%%%%%%%%%%%%%%%%%%%%%%%%%%%%%%%%%%%%%%%%
\subsubsection{Positive null controllability result}
%%%%%%%%%%%%%%%%%%%%%%%%%%%%%%%%%%%%%%%%%%%%%%%%%%%%%%%%%%%%%%%%%%%%%%%%%%%%%%%%%  
In this subsection, our objective is to prove that system  (\ref{boun.sg2x2}) is exactly controllable to zero at time $T$ if $T > T_0\in [0,\infty)$ (see (\ref{T0+})). To this end, for $y_0 \in H^{-1}(0,\pi; \mathbb{R}^2)$ we will reformulate the null controllability problem as a moment problem. Let us first observe that condition (\ref{cndca+}) is a necessary condition for having the null controllability property of system (\ref{boun.sg2x2}) at time $T > 0$. Using Proposition \ref{prop1} and Proposition \ref{prop2} we deduce that the control $u \in L^2(0, T)$ drives the solution of (\ref{boun.sg2x2})   to zero at time $T$ if and only if $u \in L^2(0, T)$ satisfies

\begin{equation*}%\label{ipp}
\int_{0}^{T}  u(t)B^{\ast }D\theta_x(0,t) dt = - \langle y^0 , \theta (\cdot , 0) \rangle_{H^{-1},H_0^1} , \quad \forall
\theta ^{0}\in H^{1}_0 ( 0,\pi ;\mathbb{R}^{2} ) ,
\end{equation*}
where $\theta$ is the solution to the adjoint problem (\ref{Padjoint}) associated with $\theta_0$. Since $\mathcal{B^{\ast}}$ is complete in $ H^{1}_0 ( 0,\pi ;\mathbb{R}^{2} )$, the null controllability property at time $T$ for system (\ref{boun.sg2x2}) is equivalent to find $u\in L^2(0, T)$ such that

\begin{equation*}
%\label{ipp+ en +}
\left\{ 
	\begin{aligned}
&\int_{0}^{T}  u(t)B^{\ast }D\theta_x^{2,j}(0,t) dt = - \langle y^0 , \theta^{2,j} (\cdot , 0) \rangle_{H^{-1},H_0^1}, &j& \in \mathbb{N}^{\ast} \setminus j_0\mathbb{N}^{\ast}, \\
&\int_{0}^{T}  u(t)B^{\ast }D\widetilde{\theta}_x^{1,k}(0,t) dt = - \langle y^0 ,\widetilde{\theta}^{1,k} (\cdot , 0) \rangle_{H^{-1},H_0^1}, &k& \in\mathbb{N}^{\ast} \setminus i_0\mathbb{N}^{\ast},
\\
&\int_{0}^{T}  u(t)B^{\ast }D\theta_x^{2,j_0l}(0,t) dt = - \langle y^0 , \theta^{2,j_0l} (\cdot , 0) \rangle_{H^{-1},H_0^1}, &l&\ge 1,\\
&\int_{0}^{T}  u(t)B^{\ast }D\widehat{\theta}_x^{1,l}(0,t) dt = - \langle y^0 , \widehat{\theta}^{1,l} (\cdot , 0) \rangle_{H^{-1},H_0^1},
 &l&\ge 1,
\end{aligned}%
\right. 
\end{equation*}
where $\theta^{2,j}$ is the solution of adjoint problem (\ref{Padjoint}) associated with $\theta^0=\Phi^{\ast}_{2,j}$, and $\widetilde{\theta}^{1,k}$ (resp. $\widehat{\theta}^{1,l}$) is the solution of adjoint problem(\ref{Padjoint}) associated with $\theta^0=\widetilde{\Phi}^{\ast}_{1,k}$  (resp. $\theta^0=\widehat{\Phi}^{\ast}_{1,k}$) . One has:
\begin{enumerate}
\item 
If we take $\theta ^{0}\equiv\Phi _{2,j}^{\ast }$, the solution of the adjoint problem is  $\theta^{2,j}(\cdot, t)= e^{-\nu j^{2}( T-t ) }\Phi _{2,j}^{\ast }$ and we obtain: 
\begin{equation*}%\label{mom1}
\int_0^T e^{-\nu j^2t} u(T-t)\, dt=e^{-\nu j^2T}\langle y^0,\frac{\Phi_{2,j}^{\ast }}{\nu \varphi'_j(0)} \rangle_{-1,1}= e^{-\nu j^2T}M(y^0, \nu j^2) , \: \forall j \in \mathbb{N}^{\ast} \setminus j_0\mathbb{N}^{\ast},
\end{equation*}
where
\begin{equation*}%\label{M}
M(y^0,\nu j^2)=\langle y^0,\frac{\Phi_{2,j}^{\ast }}{\nu \varphi'_j(0)} \rangle_{H^{-1},H^1_0}, \:\forall j \in \mathbb{N}^{\ast} \setminus j_0\mathbb{N}^{\ast}.
\end{equation*}
Moreover 

\begin{equation}
\label{M1est}
\left \vert M(y^0,\alpha) \right \vert \le C \left \Vert y_0 \right \Vert_{H^1_0(0,\pi; \mathbb{R}^2)} ,  \quad \forall \alpha \in \Lambda_1 \quad \hbox{(see (\ref{Lambdas}))},
\end{equation}
for a positive constant $C$ independent of $\alpha \in \Lambda_1$ and $y^0$.

\item 
If we take $\theta ^{0}\equiv\widetilde{\Phi} _{1,k}^{\ast }$, the solution of the adjoint problem is  $\widetilde{\theta}^{1,k}(\cdot, t)=e^{-k^2( T-t ) }\widetilde{\Phi}_{1,k}^{\ast }$ and we obtain: 
\begin{equation*}%\label{momtilde}
\int_0^T e^{-k^2t} u(T-t)\, dt=e^{- k^2T}\langle y^0,\frac{\widetilde{\Phi}_{1,k}^{\ast }}{\nu \widetilde{\psi}'_k(0)} \rangle_{-1,1}=\frac{e^{-k^2T }}{I(k^2)}\widetilde{M}(y^0,k^2), \: \forall k \in \mathbb{N}^{\ast} \setminus i_0\mathbb{N}^{\ast}. 
\end{equation*}
where 
\begin{equation*}%\label{Mtilde}
\widetilde{M}(y^0,k^2)=\nu \sqrt{\pi/2}\sin\left(\frac{k\pi}{\sqrt{\nu}}\right)\left\langle y^0,\widetilde{\Phi}_{1,k}^{\ast }\right\rangle_{H^{-1},H^1_0},  \quad \forall k \in \mathbb{N}^{\ast} \setminus i_0\mathbb{N}^{\ast}.
\end{equation*}
Using the properties of the function  $\widetilde{\psi}'_k$ stated in Proposition \ref{prop2est}, one has

\begin{equation}
\label{Mtildeest}
\vert\widetilde{M}(y^0,\beta) \vert \le C \Vert y^0\Vert_{H^1_0(0,\pi; \mathbb{R}^2)},\quad \forall \beta \in \Lambda_2 ,
\end{equation}
for a new positive constant $C$ independent of $\beta \in \Lambda_2$ and $y^0$.

\item
Let us now take $\theta ^{0}\equiv\Phi _{2,j_0l}^{\ast }$ $\left(resp. \:\theta ^{0}\equiv\widehat{\Phi}_{1,l}^{\ast }\right)$, the solution of the adjoint problem is

$\theta^{2,j_0l}(\cdot, t)=e^{-i_0^2 l^{2}( T-t ) }\Phi _{2,j_0l}^{\ast }\left(resp.\: \widehat{\theta}^{1,l}(\cdot, t)=e^{-i_0^2l^{2} ( T-t ) }\left[\widehat{\Phi}_{1,l}^{\ast }-( T-t )I(i_0^2l^2)\Phi_{2,j_0l}^{\ast }\right] \right)$, for all  $l\ge 1$. Thus, the control $u$ must also satisfy
\begin{equation*}
%\label{momhat}
\left\{ 
	\begin{aligned}
&	\int_0^T e^{-i_0^2 l^2t} u(T-t)\, dt=e^{-i_0^2 l^2 T}M^*(y^0, i_0^2l^2),\\
&\int_0^T te^{-i_0^2l^2 t} u(T-t)\, dt=\frac{e^{- i_0^2l^2T}}{I(i_0^2l^2)}\widehat{M}(y^0, i_0^2l^2), 
\end{aligned}%
\quad \forall l \ge 1,
\right. 
\end{equation*}
where %$M^{(j)}(y^0)$ and $\widehat{M}^{(l)}(y^0)$ satisfying (\ref{M1est}) and (\ref{M1lhat})
 
\begin{equation*}
%\label{M*}
M^*(y^0, i_0^2l^2)=\langle y^0,\frac{\Phi_{2,j_0l}^{\ast }}{\nu \varphi'_{j_0l}(0)} \rangle_{H^{-1},H^1_0}, \quad \forall l\ge 1,
\end{equation*}
 and, for all $l \ge 1$
\begin{equation*}
%\label{Mhat}
\widehat{M}(y^0, i_0^2l^2)=\frac{1}{\nu \varphi'_{j_0l}(0)}\left\langle y^0,\widehat{\Phi}_{1,l}^{\ast } +\left(\frac{\widehat{\psi}_l'(0)}{\varphi'_{j_0l}(0)}-T I(i_0^2l^2)\right)\Phi_{2,j_0l}^{\ast }\right\rangle_{-1,1}.
\end{equation*}
Using the properties of the function  $\widehat{\psi}'_k$ stated in Proposition \ref{prop2est} , one has
\begin{equation}
\label{M2kest}
\left\{ 
\begin{aligned}
&\vert \widehat{M}(y^0, \gamma) \vert  \le  C \Vert y^0\Vert_{H^1_0(0,\pi; \mathbb{R}^2)},\quad \forall \gamma \in \Lambda_3,\\
&\vert M^*(y^0, \gamma) \vert \le  C \Vert y^0\Vert_{H^1_0(0,\pi; \mathbb{R}^2)},\quad \forall \gamma \in \Lambda_3,
\end{aligned}
\right.
\end{equation}
for a new positive constant $C$ independent of $ \gamma \in \Lambda_3$ and $y^0$.
\end{enumerate}
Summarizing, we have proved that $u\in L^2(0; T)$ is such that the solution $y$ of system (\ref{boun.sg2x2}) satisfies $y(\cdot , T) = 0$ in $(0; \pi)$ if and only if, $u\in L^2(0; T)$ satifies

\begin{equation}
\label{summarizing}
\left\{ 
	\begin{aligned}
&\int_0^T e^{-\alpha t} u(T-t)\, dt=e^{-\alpha T}M(y^0, \alpha ) , &&\forall\alpha \in \Lambda_1,\\
&\int_0^T e^{-\beta t} u(T-t)\, dt=\frac{e^{-\beta T }}{I(\beta)}\widetilde{M}(y^0,\beta), && \forall \beta\in \Lambda_2,\\
&	\int_0^T e^{-\gamma t} u(T-t)\, dt=e^{-\gamma T}M^*(y^0,\gamma), &&\forall \gamma \in \Lambda_3,\\
&\int_0^T te^{-\gamma t} u(T-t)\, dt=\frac{e^{- \gamma T}}{I(\gamma)}\widehat{M}(y^0,\gamma),
 && \forall \gamma \in \Lambda_3.
\end{aligned}%
\right.
\end{equation}

\begin{prpstn}[gap-condition]
\label{gapp}
\begin{equation*}
% \label{gap}
\forall \, \alpha,\beta \, \in \sigma(L^{\ast}),\quad \hbox{with}\quad \alpha \neq \beta,\quad \left \vert  \alpha- \beta \right \vert > \frac{1}{j_0^2}.
\end{equation*}
\end{prpstn}
\noindent
\textit {Proof of Proposition \ref{gapp}}. 
 For all $ \alpha,\beta \, \in \sigma(L^{\ast})$, there exist $p, q\in \mathbb{N}^{\ast}$ where $p \neq q$, such that $j_0^2\alpha=p^2$ and $j_0^2\beta=q^2$ (see (\ref{Lambdas}) and (\ref{sigmaLast})), consequently

$$
\left \vert   \alpha -  \beta \right \vert = \frac{1}{j_0^2}\left \vert  p^2-q^2 \right \vert \ge  \frac{1}{j_0^2}.
$$
$\hfill \square$

Thus, from the results in \cite{FGT} (see Lemma 3.1. (b), page 1725-1726), we know that the sequence
$$
\{ e^{-\alpha t},te^{-\alpha t} \}_{\alpha \in \sigma(L^{\ast})}
$$
admits a biorthogonal family $\{ q_{1,\alpha},q_{2,\alpha}\}_{\alpha\in \sigma(L^{\ast})}$ in $L^2(0,T)$, i.e., a family satisfying

\begin{equation*}
 %\label{biorth en +}
 \left\{ 
 \begin{aligned}
  &\int_0^T e^{-\alpha t}q_{1,\beta}(t)dt=\delta_{\alpha \beta}, \\
  &  \int_0^T te^{-\alpha t}q_{1,\beta}(t)dt=0, 
    \end{aligned}
 \right.
 \quad
 \hbox{and}
 \quad
  \left\{ 
 \begin{aligned}
   &\int_0^T e^{-\alpha t}q_{2,\beta}(t)dt =0, \\
  &  \int_0^T te^{-\alpha t}q_{2,\beta}(t)dt= \delta_{\alpha \beta}, 
    \end{aligned}
 \right.
\end{equation*}
where  $\alpha,\beta \in \sigma(L^{\ast})$, which moreover satisfies that for every $\varepsilon>0$ there exists a constant $C_{\varepsilon,T}>0$ such that 
\begin{equation}
\label{biorthnur1}
\Vert q_{j,\beta} \Vert _{L^2(0,T)} \le C_{\varepsilon,T}e^{\varepsilon\beta},\quad \beta\in \sigma(L^{\ast}),\quad j\in \{ 1,2\}.
\end{equation}
Using now the formulas in (\ref{summarizing}) and the property (\ref{biorthnur1}), we infer that an explicit formal solution $u$  is given by

\begin{equation}
 \label{u}
\begin{aligned}
 u(T-t)
 &=\displaystyle{\sum_{\alpha \in \Lambda_1}} e^{-\alpha T}M(y^0,\alpha)q_{1,\alpha}(t)+ \displaystyle{\sum_{\beta \in \Lambda_2}} \frac{e^{-\beta T}\widetilde{M}(y^0,\beta)q_{1,\beta}(t)}{I(\beta)}\\
& + \displaystyle{\sum_{\gamma \in \Lambda_3}}e^{-\gamma T}\left( M^*(y^0,\gamma)q_{1,\gamma}(t)+  \frac{\widehat{M}(y^0,\gamma)q_{2,\gamma}(t)}{I(\gamma)}\right)
\end{aligned}
\end{equation}

Let us see that this series defines an element of $L^2(0, T)$ when $T >  T_{0}$ , i.e., the previous series converges in $L^2(0, T)$ if $T >  T_{0}$. Indeed, from the definition of the minimal time $T_{0}$ (see (\ref{T0alpha})) and for any fixed $\varepsilon>0$, we can infer that there exists a positive constant $C_{\varepsilon}$ such that
$$
\frac{1}{\left \vert I(\zeta) \right \vert}\le C_{\varepsilon}e^{\zeta(T_{0}+\varepsilon)} , \quad \forall\zeta \in \Lambda_2\cup \Lambda_3.
$$
On one hand, we can use (\ref{biorthnur1}), (\ref{M1est}) and (\ref{M2kest}), and get for any $\alpha \in \Lambda_1$ and $\gamma \in \Lambda_3$ , a  positive constant $C_{\varepsilon,T}>0$ for which
 $$
 %\begin{aligned}
\left \Vert  e^{-\alpha T}M(y^0,\alpha)q_{1,\alpha}(t) \right \Vert_{L^{2}(0,T)}  +\left \Vert e^{-\gamma T} M^*(y^0,\gamma)q_{1,\gamma}(t)\right\Vert_{L^{2}(0,T)} \le 
C_{\varepsilon,T}\left \Vert y^0\right\Vert_{H^{1}_0}  \left( e^{\alpha (\varepsilon-T)} +e^{\gamma(\varepsilon- T)} \right). 
%\end{aligned}
$$ 
On the other hand, we can use (\ref{biorthnur1}), (\ref{Mtildeest}) and (\ref{M2kest}), and get for any $\beta \in \Lambda_2$ and $\gamma \in \Lambda_3$ , a new positive constant $C_{\varepsilon,T}>0$ such that
 $$
\left \Vert \frac{e^{-\beta T}}{I(\beta)}\widetilde{M}(y^0,\beta)q_{2,\beta}(t) \right\Vert_{L^{2}}
\le 
C_{\varepsilon,T}\left \Vert y^0\right\Vert_{H^{1}_0}  \frac{e^{\beta(\varepsilon- T)}}{\left \vert I(\beta) \right \vert}  \le
C_{\varepsilon,T}\left \Vert y^0\right\Vert_{H^{1}_0}    e^{\beta(T_0-T+2\varepsilon)},
$$

$$
\left \Vert  \frac{e^{-\gamma T}}{I(\gamma)}\widehat{M}(y^0,\gamma)q_{2,\gamma}(t) \right\Vert_{L^{2}}
\le 
C_{\varepsilon,T}\left \Vert y^0\right\Vert_{H^{1}_0}  \frac{e^{\gamma (\varepsilon- T)}}{\left \vert I(\gamma)\right \vert} \le
C_{\varepsilon,T}\left \Vert y^0\right\Vert_{H^{1}_0}    e^{\gamma(T_0-T+2\varepsilon)}.
$$
These previous two inequalities prove the absolute convergence of the series (\ref{u}) which defines the control $u$ since $\varepsilon$ may be chosen arbitrarily small (take $\varepsilon\in (0,\frac{T-T_0}{2})$ for instance). This proves the null controllability of system~(\ref{boun.sg2x2}) at time $T$ when  $T>T_{0}$.

%%%%%%%%%%%%%%%%%%%%%%%%%%%%%%%%%%%%%%%%%%%%%%%%%%%%%%%%%%%
\subsubsection{The negative null controllability result}
%%%%%%%%%%%%%%%%%%%%%%%%%%%%%%%%%%%%%%%%%%%%%%%%%

Let us prove that if $0<T<T_0 $, then system (\ref{boun.sg2x2}) is not null controllable at time $T$. We argue by contradiction.
Let us first remark that 

\begin{equation*}
%\label{maxt1t2 en+}
T_0=max\left\lbrace T_{0,1}, T_{0,2}\right\rbrace
\end{equation*}
with

\begin{equation}
\label{t1t20}
T_{0,1}=\limsup_{k \in \mathbb{N}^{\ast}\setminus i_0\mathbb{N}^{\ast},\;k \to +\infty}\frac{\log\frac{1}{\left \vert I(k^2)\right \vert}}{k^2} \quad \hbox{and} \quad
T_{0,2}=\limsup_{l \in  \mathbb{N}^{\ast},\; l\to +\infty } \frac{\log\frac{1}{\left \vert I(i_0^2 l^2) \right \vert}}{i_0^2 l^2} .
\end{equation}
Assume  that ~\eqref{boun.sg2x2}  is  null-controllable at time $T<T_0 $,  system~(\ref{boun.sg2x2}) is null-controllable at time $T$ if and only if there exists $C>0$ such that any solution $\theta $ of the adjoint problem~(\ref{Padjoint}) satisfies the observability inequality: 
\begin{equation}
\Vert \theta (0)\Vert _{H^1_0(0,\pi ;\mathbb{R}^{2})}^{2}\leq
C\int_{0}^{T}\!\!|B^*D\theta _{x}(0,t)|^{2}\,dt,\quad \forall \theta
^{0}\in H^1_0\left( 0,\pi ;\mathbb{R}^{2}\right) . 
\label{IO}
\end{equation}
\begin{itemize}
\item  
Let us suppose   $T_0=T_{0,1}$. Let us work with the particular solutions $\theta_k$ associated with initial data $$\theta ^{0}_k=\frac{\widetilde{\Phi}_{1,k}^{\ast }}{\nu \widetilde{\psi}'_k(0)} .$$
 With this choice, the solution $\theta_k$ of (\ref{Padjoint}) is given by

$$ \theta_{k} (x,t)=\frac{\widetilde{\Phi}_{1,k}^{\ast }}{\nu \widetilde{\psi}'_k(0)} e^{-k^2(T-t)},\quad  k \in \mathbb{N}^{\ast}\setminus i_0\mathbb{N}^{\ast}.
$$ 
From the definition of $T_{0,1}$ (see (\ref{t1t20})), then there exists an  sub-sequence $\left(k_n\right)_{n \ge 1} \subset  \mathbb{N}^{\ast}\setminus i_0\mathbb{N}^{\ast}$ such that:
$$
T_{0,1}=\lim_{n\to+\infty}\frac{\log\frac{1}{\vert I(k_n^2)\vert}}{k_n^2}\in(0,+\infty],
$$ 
thus
\begin{equation*}
% \begin{cases}
 \begin{aligned}
\Vert \theta_{k_n} (\cdot,0)\Vert_{H^1_0(0,\pi ;\mathbb{R}^{2})}^2 &=
 \frac{\left \Vert \widetilde{\Phi}_{1,k_n}^{\ast } \right\Vert_{H^1_0}^2}{\left \vert \nu \widetilde{\psi}'_{k_n}(0) \right \vert^2} e^{-2k^2_n T}  = \frac{k_n^2+\left \Vert \widetilde{\psi}'_{k_n}\right\Vert_{L^2}^2}{\left \vert \nu \widetilde{\psi}'_{k_n}(0) \right \vert^2} e^{-2k^2_n T},\\
&\ge  \frac{\pi}{2}  \frac{\left \vert\sin\left( \frac{k_n \pi}{\sqrt{\nu}}\right)\right \vert^2}{\left \vert I(k_n^2) \right \vert} e^{-2k^2_n T} \ge k_n^2 \frac{\pi}{2} \frac{4}{i_0^2}  \frac{1}{\left \vert I(k_n^2) \right \vert^2}e^{-2k^2_n T} ,\\
&\ge  k_n^2 \frac{\pi}{2} \frac{4}{i_0^2}e^{-2k^2_n  \left[ T-\frac{ \log\left \vert I(k_n^2) \right \vert }{ k_n^2} \right]} ,  \\
&\ge  k_n^2 \frac{\pi}{2} \frac{4}{i_0^2}e^{-2k^2_n \left[ T-T_{0,1} +\varepsilon\right]} \underset{n\rightarrow +\infty}{\longrightarrow} +\infty, \quad \varepsilon\in(0,T_{0,1}-T),
\end{aligned}
%\label{NI1}
\end{equation*}
and
\begin{equation*}
\int_{0}^{T}|B^{\ast }D \frac{\partial \theta_{k_n} (0,t)}{\partial x} (0,t)|^{2}\,dt=\int_{0}^{T} e^{-2k_n^2t}\,dt \underset{n\rightarrow +\infty}{\longrightarrow} 0.
%\end{cases}
% \label{NI1+}
\end{equation*}
This proves that the observability inequality (\ref{IO}) does not hold.

\item Let us suppose now, $T_0=T_{0,2}$ . Let us work with the particular solutions $\theta^l$ associated with initial  data 
$$\theta ^{0}=a_{l}\widehat{\Phi} _{1,l}^{\ast }+b_{l}\Phi _{2,j_0l}^{\ast },\quad l \ge 1,$$
  with  $a_{l},b_{l} \in \mathbb{R}$, to be determined. Recall that  
  $$\int_{0}^{\pi}\widehat{\psi}_l(x) \varphi_{j_0l}(x)\,dx=0\quad \hbox{(see (\ref{psikeq}))}.$$
With this choice, the solution $\theta^l$ of (\ref{Padjoint}) is given by

$$
\theta^{l}(\cdot, t)=a_le^{-i_0^2 l^2 ( T-t ) }\left[\widehat{\Phi} _{1,l}+( T-t )I(i_0^2 l^2)
\Phi_{2,j_0l}^{\ast }\right]+b_le^{-i_0^2 l^2( T-t ) }\Phi_{2,j_0l}^{\ast },\quad \forall l\ge 1,$$
thus, the observability inequality (\ref{IO})  becomes
$$
A_{1,l}\le CA_{2,l}
$$
with, for all $l\ge 1$
$$
 \begin{cases}
 \begin{aligned}
& A_{1,l}=e^{-2i_0^2 l^2T}\left\{i_0^2 l^2|a_{l}|^{2}+\left[|a_{l}|^{2}\Vert\psi'_{l}\Vert_{L^{2}(0,\pi)}^{2}+i_0^2 l^2\left(b_{l}+Ta_{l}I(i_0^2 l^2)\right)^{2}\right]\right\}\ge  e^{-2i_0^2 l^2T}|a_{l}|^{2},  \\
 & A_{2,l}=\int_{0}^{T}e^{-2i_0^2 l^2t}\left\vert a_{l}\psi'_{l}(0)+(b_{l}+ta_{l}I(i_0^2 l^2))\varphi'_{j_0l}(0)\right\vert ^{2}\,dt.
 \end{aligned}
\end{cases}.
$$
We can choose $a_{l}=I(i_0^2 l^2)$ and $b_{l}=-I(i_0^2 l^2)\psi'_{l}(0)/\varphi'_{j_0l}(0)$, we obtain for a new constant $C > 0$ not depending on $l$, such that: 

$$
\frac{A_{1,l}}{A_{2,l}}\ge \frac{ e^{-2i_0^2 l^2T}}{\left\vert I(i_0^2 l^2)\varphi'_{j_0 l}(0)\right\vert ^{2}\int_{0}^{T}t^2e^{-2i_0^2 l^2t}\,dt} \ge \frac{ e^{-2i_0^2 l^2 T}}{\left\vert I(i_0^2 l^2)\right \vert^2},\quad \forall l\ge 1.
$$
From the definition of $T_{0,2}$ (see (\ref{t1t20})), there exists a sub-sequence  $\{l_{n}\}_{n\ge1}$
 such that : 
$$
T_{0,2}=\lim_{n\to+\infty}\frac{\log\frac{1}{\vert I(i_0^2 l_n^2)\vert}}{i_0^2 l_n^2}\in(0,+\infty].
$$ 
Assume $0<T_{0,2}<+\infty$. In this case, for every $\varepsilon>0$, there exits a positive integer $n_{\varepsilon}\geq1$ such that 

$$ 
 \frac{1}{\vert I(i_0^2 l_n^2)\vert}>e^{i_0^2 l_n^2(T_0^2-\varepsilon)},\quad\forall n\ge n_{\varepsilon},
 $$
thus, $\forall n\ge n_{\varepsilon}$
 
$$
\frac{A_{1,l_n}}{A_{2,l_n}} \ge \frac{ e^{-2i_0^2 l_n^2 T}}{\left\vert I(i_0^2 l_n^2)\right \vert^2}\ge  e^{2i_0^2 l_n^2(T_{0,2}-T-\varepsilon)}.
$$
Taking for instance $\varepsilon=\frac{T_{0,2}-T}{2}$, we have
$$
e^{2i_0^2 l_n^2(T_{0,2}-T-\varepsilon)}=e^{i_0^2 l_n^2(T_{0,2}-T)}\underset{n\rightarrow +\infty}{\longrightarrow}+\infty.
 $$

\end{itemize}
This proves that the observability inequality (\ref{IO}) does not hold and finishes the proof of the negative null controllability result of (\ref{boun.sg2x2}).

%%-----------------------------
%%      your bibliography
%%-----------------------------

%%%%%%%%%%%Appendix
\appendix
\section{Reformulation of minimal time in Theorem \ref{Thm} }
\label{AppA}
%%%%%%%%%%%
%%%%%%%%%%%%%%%%%%%%%%%%%%%%%%%%%%%%%%%%%%%%%%%%%%%%%%%%%%%%%%%%%%%
%\section{ Reformulations des temps minimaux du Théorème \ref{Thm}}
%%%%%%%%%%%%%%%%%%%%%%
Let us begin this section by reformulating the minimal times introduced in Theorem \ref{Thm}. 
Let us first recall that 
\begin{equation} 
\label{app1}
B^{\ast}D\frac{\partial\Phi_{1,k}^{\ast }}{\partial x} (0)=\psi'_{k}(0)=-\frac{\int_{0}^{\pi}q(s)\varphi_{k}(s)\sin\left(\frac{k}{\sqrt{\nu}}(\pi-s)\right)ds}{\nu \sin\left(\frac{k \pi}{\sqrt{\nu}}\right)}.
\end{equation}
In all the sequel, assume that condition (\ref{cndca}) is holds. We have

\begin{prpstn}
Let us consider the definitions of the Theorem \ref{Thm}. One has
\label{propo reform tm}
\begin{enumerate}
    \item
    \begin{equation*}
        T_1 =\displaystyle{\limsup_{k \to +\infty } } \left\{\frac{1}{k^2}\log \frac{\left \vert  \sin\left(\frac{k \pi}{\sqrt{\nu}}\right)\right \vert}{\left \vert \int_{0}^{\pi}q(s)\varphi_{k}(s)\sin\left(\frac{k}{\sqrt{\nu}}(\pi-s)\right)ds\right \vert}\right\}\, \in [0, +\infty].
    \end{equation*}
    
 If  $\sqrt{\nu}>1$, then
      \begin{equation*}
        \widetilde{T}_2= \displaystyle{\limsup_{k \to +\infty } } \left\{ \frac{1}{i_k^2}\log\frac{1}{ \left \vert\int_{0}^{\pi}q(s)\varphi_{i_k}(s)\sin\left(\frac{i_k}{\sqrt{\nu}}(\pi-s)\right)ds \right\vert} \right\}\, \in [0, +\infty].
    \end{equation*}

      \item If  $\sqrt{\nu}<1$, then
   
      \begin{equation*}
        \widehat{T}_2= \displaystyle{\limsup_{k \to +\infty } }\left\{ \frac{1}{k^2}\log\frac{1}{ \left \vert  \int_{0}^{\pi}q(s)\varphi_{k}(s)\sin\left(\frac{k}{\sqrt{\nu}}(\pi-s)\right)ds\right  \vert}\right\}\ge T_1 , 
    \end{equation*}
    consequently
 \begin{equation*}
    \widehat{T}_0=\widehat{T}_2.
     \end{equation*}
    
\end{enumerate}
\end{prpstn}

\bigskip

\noindent
\textit {Poof of Proposition \ref{propo reform tm} }.
Let us recall that for all $k\ge 1$, one has
\begin{equation*} \Phi _{1,k}^{\ast }:=\left( 
\begin{array}{l}
\varphi _{k} \\ 
\psi _{k}%
\end{array}%
\right) ,\: \Phi _{2,k}^{\ast }:=\left( 
\begin{array}{l}
0 \\ 
\varphi _{k}%
\end{array}%
\right) \quad \hbox{(see (\ref{vectpropnuir}))},
\end{equation*}%
and
\begin{equation*} 
\psi_{1,k}:=\frac{\Phi_{1,k}^{\ast }}{B^{\ast}D\frac{\partial\Phi_{1,k}^{\ast }}{\partial x} (0)}=\frac{\Phi_{1,k}^{\ast }}{\nu \psi'_{k}(0)} \quad \hbox{(see (\ref{psik}))},
\end{equation*}
where (see (\ref{tpsiksolnuir}))
\begin{equation*} 
\psi_{k}(x)=\psi'_{k}(0)\frac{\sqrt{\nu}}{k}\sin\left(\frac{kx}{\sqrt{\nu}}\right)+\frac{\sqrt{\nu}}{\nu k}\int_{0}^{x}\sin\left(\frac{k}{\sqrt{\nu}}(x-\xi)\right)q(\xi)\varphi_{k}(\xi)\,d\xi.
\end{equation*}%
Moreover, one has
\begin{equation*} 
\psi'_{k}(x)=\psi'_{k}(0)\cos\left(\frac{kx}{\sqrt{\nu}}\right)+\frac{1}{\nu }\int_{0}^{x}\cos\left(\frac{k}{\sqrt{\nu}}(x-\xi)\right)q(\xi)\varphi_{k}(\xi)\,d\xi.
\end{equation*}%
We immediately see that there is a constant $C>0$ such that 

\begin{equation}
\Vert \psi'_k\Vert^2_{L^2(0,\pi)} \le C \vert \nu \psi'_k (0)\vert^2+C.
\label{est_psik_exp}
\end{equation}%

\begin{enumerate}
    \item 
    On one hand, using the previous inequality (\ref{est_psik_exp}), one has
    $$
    \begin{aligned}
        \frac{k^2}{\vert \nu \psi'_k(0)\vert^2} & \le  \Vert \psi_{1,k} \Vert^2_{H^1_0} =\frac{k^2+ \Vert \psi'_k \Vert^2_{L^2}}{\vert \nu \psi'_k(0)\vert^2}\le  \frac{k^2+ C \vert \nu \psi'_k (0)\vert^2+C}{\vert \nu \psi'_k(0)\vert^2}\le  \frac{k^2+ C}{\vert \nu \psi'_k(0)\vert^2}+C,
    \end{aligned}
    $$
   % Comme
  %  \begin{equation*}
   % \label{ineg sinus}
  %2\left \vert \frac{k }{\sqrt{\nu}} -j_k \right\vert \le  \left \vert  \sin\left(\frac{k \pi}{\sqrt{\nu}}\right)\right \vert \le \pi \left \vert \frac{k }{\sqrt{\nu}} -j_k \right\vert,
  %  \end{equation*}%
   %  où $j_k$ est l'entier le plus proche de $\frac{k}{\sqrt{\nu}}$, alors
   thus 
   $$
   \begin{aligned}
     T_1 &=\displaystyle{\limsup_{k \to +\infty } } \frac{\log \left \Vert \psi_{1,k}\right \Vert_{H^1_0}}{k^2}=\displaystyle{\limsup_{k \to +\infty } } \frac{-\log \left \vert \nu \psi'_{k}(0)\right \vert}{k^2}\\
     &=\displaystyle{\limsup_{k \to +\infty } } \left\{\frac{1}{k^2}\log \frac{\left \vert  \sin\left(\frac{k \pi}{\sqrt{\nu}}\right)\right \vert}{\left \vert \int_{0}^{\pi}q(s)\varphi_{k}(s)\sin\left(\frac{k}{\sqrt{\nu}}(\pi-s)\right)ds\right \vert}\right\} \quad \hbox{(see (\ref{app1}))}.
      \end{aligned}
    $$
 On the other hand, for all $k\ge 1$ and $\sqrt{\nu}>1$, one has
    
    \begin{equation}
   \label{inepp1}
    \vert i_k-\sqrt{\nu}k \vert <\frac{1}{2},
     \end{equation}
    where $i_k$  is the nearest integer to $\sqrt{\nu}k$, thus
    $\vert \frac{i_k}{\sqrt{\nu}}-k \vert <\frac{1}{2\sqrt{\nu}}<\frac{1}{2}$. Consequently
    \begin{equation}
    \label{ineg sinus}
  2\left \vert \frac{i_k }{\sqrt{\nu}} -k \right\vert \le  \left \vert  \sin\left(\frac{i_k \pi}{\sqrt{\nu}}\right)\right \vert \le \pi \left \vert \frac{i_k }{\sqrt{\nu}} -k \right\vert.
    \end{equation}
Thus
  
   $$
    \begin{aligned}
   \frac{\left \Vert \psi_{1,i_k}-\psi_{2,k} \right \Vert^2_{H^1_0}}{\vert  i_k^2-\nu k^2\vert^2 } & =   \frac{1}{\vert  i_k^2-\nu k^2\vert^2 } \left(\frac{i_k^2}{\vert \nu \psi'_{i_k}(0)\vert^2} + \frac{ \Vert \nu \varphi'_k(0)\psi'_{i_k} - \nu \psi'_{i_k}(0)\varphi'_{k}\Vert^2_{L^2}}{\vert \nu \psi'_{i_k}(0)\vert^2\vert \nu \varphi'_k(0)\vert^2} \right)\\
      & \le  \frac{1}{\vert  i_k^2-\nu k^2\vert^2 } \left(\frac{i_k^2}{\vert \nu \psi'_{i_k}(0)\vert^2} + \frac{2\Vert \psi'_{i_k}\Vert^2_{L^2}}{\vert \nu \psi'_{i_k}(0)\vert^2} +\frac{2k^2}{\vert \nu \varphi'_k(0)\vert^2} \right)\\
      & \le \frac{1}{\vert  i_k^2-\nu k^2\vert^2 } \left(\frac{i_k^2+C}{\vert \nu \psi'_{i_k}(0)\vert^2} +C +\frac{\pi}{ \nu } \right)\\
      & \le  \frac{C'i_k^2}{\vert  i_k^2-\nu k^2\vert^2 } \left(\frac{ \left \vert  \sin\left(\frac{i_k \pi}{\sqrt{\nu}}\right)\right \vert^2}{\left \vert\int_{0}^{\pi}q(s)\varphi_{i_k}(s)\sin\left(\frac{i_k}{\sqrt{\nu}}(\pi-s)\right)ds \right\vert^2}  \right)\\
      & \le  \frac{C'i_k^2}{\vert  i_k+\sqrt{\nu} k\vert^2 }\frac{\pi^2}{\nu} \left(\frac{1}{\left \vert\int_{0}^{\pi}q(s)\varphi_{i_k}(s)\sin\left(\frac{i_k}{\sqrt{\nu}}(\pi-s)\right)ds \right\vert^2}  \right),
    \end{aligned}
    $$
and
 
   $$
    \begin{aligned}
   \frac{\left \Vert \psi_{1,i_k}-\psi_{2,k} \right \Vert^2_{H^1_0}}{\vert  i_k^2-\nu k^2\vert^2 } & =   \frac{1}{\vert  i_k^2-\nu k^2\vert^2 } \left(\frac{i_k^2}{\vert \nu \psi'_{i_k}(0)\vert^2} + \frac{ \Vert \nu \varphi'_k(0)\psi'_{i_k} - \nu \psi'_{i_k}(0)\varphi'_{k}\Vert^2_{L^2}}{\vert \nu \psi'_{i_k}(0)\vert^2\vert \nu \varphi'_k(0)\vert^2} \right)\\
      & \ge  \frac{1}{\vert  i_k^2-\nu k^2\vert^2 } \frac{i_k^2}{\vert \nu \psi'_{i_k}(0)\vert^2} \\
      & \ge  \frac{i_k^2}{\vert  i_k^2-\nu k^2\vert^2 } \left(\frac{ \left \vert  \sin\left(\frac{i_k \pi}{\sqrt{\nu}}\right)\right \vert^2}{\left \vert\int_{0}^{\pi}q(s)\varphi_{i_k}(s)\sin\left(\frac{i_k}{\sqrt{\nu}}(\pi-s)\right)ds \right\vert^2}  \right)\\
      & \ge \frac{i_k^2}{\vert  i_k+\sqrt{\nu} k\vert^2 } \frac{4}{\nu}\left(\frac{1}{\left \vert\int_{0}^{\pi}q(s)\varphi_{i_k}(s)\sin\left(\frac{i_k}{\sqrt{\nu}}(\pi-s)\right)ds \right\vert^2}  \right),
    \end{aligned}
    $$
and then, we deduce that
    $$
     \widetilde{T}_2= \displaystyle{\limsup_{k \to +\infty } } \frac{-\log \left \vert\int_{0}^{\pi}q(s)\varphi_{i_k}(s)\sin\left(\frac{i_k}{\sqrt{\nu}}(\pi-s)\right)ds \right\vert}{\nu k^2}.
    $$ 
  
  \item
 Similarly, we show that,  for $\sqrt{\nu}<1$

    $$
     \widehat{T}_0= \widehat{T}_2= \displaystyle{\limsup_{k \to +\infty } } \frac{-\log \left \vert\int_{0}^{\pi}q(s)\varphi_{k}(s)\sin\left(\frac{k}{\sqrt{\nu}}(\pi-s)\right)ds \right\vert}{k^2} \ge T_1.
    $$

\end{enumerate}
 $\hfill \square $
%%%%%%%%%%
%\section{Quelques exemples}
%%%%%%%%%%%%%%%%%%%%%%

In the following Proposition \ref{prop exmple} we will give a positive answer to the natural question : Let $\tau_0 \in[0,+\infty]$, is there $\sqrt{\nu} \notin \mathbb{Q}^*$ and $q\in L^{\infty}(0,\pi)$ such that $\widetilde{T}_0=\tau_0$ (or  $\widehat{T}_0=\tau_0$) ?

\begin{prpstn}
 Let $\sqrt{\nu }\notin \mathbb{Q}$.  
   \begin{enumerate}
    \item For any $\tau_0 \in[0,+\infty]$, there exists $\sqrt{\nu} >1$ and  $q\in L^{\infty}(0,\pi)$ such that
    $$
    \widetilde{T}_0=\max\{ T_1,\widetilde{T}_2 \}=\tau_0.
    $$
 \item For any $\tau_0 \in[0,+\infty]$, there exists $\sqrt{\nu}<1$ and  $q\in L^{\infty}(0,\pi)$ such that
    $$
    \widehat{T}_0=\widehat{T}_2=\tau_0.
    $$
   
\end{enumerate}
\label{prop exmple}
\end{prpstn}

\noindent
\textit {Proof of Proposition \ref{prop exmple} }
%\begin{enumerate}
  %  \item 
In the sequel, we will work with $\widetilde{T}_0$. The arguments to prove that $\widehat{T}_0=\tau_0\in [0,+\infty]$, are similar. Let us recall

 \begin{lmm}
\begin{enumerate}
\item Let us fixed $\tau_0 \in (0, \infty)$, $x_0 \in [0, \infty)$ and $%
\varepsilon > 0$. Then, there exist an irrational number  $\nu >0$  and a sequence of rational numbers
  $\left\{ p_{k}/q_{k}\right\}_{k \geq 1}$ such that $p_k$ and $q_k$ are co-prime positive integers, the sequences  $\left\{p_{k}\right\}_{k \geq 1}$ and $\left\{ q_{k}\right\}_{k \geq 1}$ are strictly increasing,, 
\begin{equation}  
\label{eq1lem1}
| \nu - x_0 | \le \varepsilon \quad \hbox {and} \quad \lim_{k \to \infty} e^{\tau_0
q_{k}^{2}}\left\vert \nu -\frac{ p_{k}}{q_{k}}\right\vert =1.
\end{equation}
Moreover, for any $k \ge 1$  one has
\begin{equation}  
\label{eq2lem1}
0 < \left| q_k \nu - p_k \right| \le \left| q \nu - p \right| , \quad
\forall p,q \in \mathbb{N}^*, \quad \hbox{with } q < q_{k+1}.
\end{equation}

\item  For any  $\sigma \in (0, \infty)$, $x_0 \in [0, \infty)$ and $%
\varepsilon > 0$, there exists an irrational number  $\nu >0$  and a
sequence of rational numbers $\left\{ p_{k}/q_{k}\right\}_{k \geq 1}$ such that  $p_k$ and $q_k$ are co-prime positive integers, the sequences $\left\{p_{k}\right\}_{k \geq 1}$ and $\left\{ q_{k}\right\}_{k \geq 1}$ are strictly increasing and 
\begin{equation}  \label{eq3lem1}
| \nu - x_0 | \le \varepsilon \quad \hbox {and} \quad  \lim_{k \to \infty} 
e^{q_{k}^{2+\sigma }}\left\vert \nu - \frac{p_{k}}{q_{k}}\right\vert =0.
\end{equation}
\end{enumerate}
 \label{lem1}
\end{lmm}

\noindent

\textit {Proof of Lemma \ref{lem1} } See for example  \cite{JMAA16} Lemma A.1, page 35. 
$\hfill \square $

Note that by choosing $q(x)=1,\: x \in [0,\pi]$, we obtain
$$
%\begin{array}{ccc}
    \frac{\left \vert  \sin\left(\frac{k \pi}{\sqrt{\nu}}\right)\right \vert}{\left \vert \int_{0}^{\pi}q(s)\varphi_{k}(s)\sin\left(\frac{k}{\sqrt{\nu}}(\pi-s)\right)ds\right \vert}  =   \frac{\left \vert  \sin\left(\frac{k \pi}{\sqrt{\nu}}\right)\right \vert}{\left \vert \int_{0}^{\pi}\varphi_{k}(s)\sin\left(\frac{k}{\sqrt{\nu}}(\pi-s)\right)ds\right \vert}=\sqrt{\frac{\pi}{2}}\frac{k\left\vert \nu-1\right\vert}{\nu},
%\end{array}
$$
and thus $T_1=0$. But 

$$
%\begin{array}{ccc}
    \frac{1}{\left \vert \int_{0}^{\pi}q(s)\varphi_{i_k}(s)\sin\left(\frac{i_k}{\sqrt{\nu}}(\pi-s)\right)ds\right \vert}  =   \frac{1}{\left \vert \int_{0}^{\pi}\varphi_{i_k}(s)\sin\left(\frac{i_k}{\sqrt{\nu}}(\pi-s)\right)ds\right \vert}=\sqrt{\frac{\pi}{2}}\frac{i_k\vert \nu-1\vert}{\nu\left \vert  \sin\left(\frac{i_k \pi}{\sqrt{\nu}}\right)\right \vert},
%\end{array}
$$
thus in this case
$$
\widetilde{T}_0=\widetilde{T}_2= \displaystyle{\limsup_{k \to +\infty } } \frac{-\log \left \vert\sin\left(\frac{i_k\pi}{\sqrt{\nu}}\right)\right\vert}{i_k^2}.
$$
Finally, we have
\begin{lmm}
\label{egalitésinus}
Let $\nu \notin\mathbb{Q}^*$ and $i_k$ nearest integer to $\sqrt{\nu}k$, for each $k\ge 1$. If $\nu>1$ then
$$
\displaystyle{\limsup_{k \to +\infty } } \frac{-\log \left \vert\sin\left(\frac{i_k\pi}{\sqrt{\nu}}\right)\right\vert}{i_k^2}=\displaystyle{\limsup_{k \to +\infty } } \frac{-\log \left \vert\sin\left(\sqrt{\nu}k \pi \right)\right\vert}{\nu k^2}.
$$
\end{lmm}

\noindent
\textit {Proof of Lemma \ref{egalitésinus} }
We have

$$
\begin{aligned}
  \displaystyle{\limsup_{k \to +\infty } } \frac{-\log \left \vert\sin\left(\frac{i_k\pi}{\sqrt{\nu}}\right)\right\vert}{i_k^2} & =  \displaystyle{\limsup_{k \to +\infty } } \frac{-\log \left \vert \frac{i_k}{\sqrt{\nu}}- k\right\vert}{\nu k^2},\quad \left( \hbox{see}\: (\ref{ineg sinus}), \: \hbox{and} \: \lim \frac{i_k^2}{\nu k^2}=1 \right) \\
   & =  \displaystyle{\limsup_{k \to +\infty } } \frac{-\log \left \vert \sqrt{\nu}k-i_k\right\vert}{\nu k^2}, \\
   & =  \displaystyle{\limsup_{k \to +\infty } } \frac{-\log \left \vert \sin(\sqrt{\nu}k)\right\vert}{\nu k^2}\quad  \left(\hbox{see}\: (\ref{inepp1}) \right).
\end{aligned}
$$
 $\hfill \square $

Let us fixed $x_0 \ge 0$ et $\varepsilon>0$. We will distinguish three different cases in the proof of Proposition \ref{prop exmple}. 
\textbf{Case $\tau_0=0$ :} 
Given $x_0 \ge 1$ and $\varepsilon > 0$ such that  $x_0-\varepsilon >1$.
Taking  $\nu>1$  a irrational algebraic number of order $n\ge 2$, we deduce the existence of a positive constant $C>0$ such that 
$$
\vert \sin(i_k\pi/\sqrt{\nu} \vert \ge \frac{C}{i_k^{n-1}}, \quad \forall k\ge 1,
$$
Consequently $\widetilde{T}_0=\widetilde{T}_2=0$.

\textbf{Case $\tau_0\in (0,\infty)$ :}
Given $\tau_0$, $x_0\ge 0$ and $\varepsilon > 0$ such that  $x_0-\varepsilon >1$. We deduce the existence of a irrational algebraic number  $\nu >1$  which fulfills (\ref{eq1lem1}) and
(\ref{eq2lem1}) for the sequences of positive integers $\{p_k\}_{k \ge 0}$ et $\{q_k \}_{k\ge 0}$. From (\ref{eq1lem1}) we deduce $\lim(p_k-q_k\sqrt{\nu})=0$, $\lim p_k/q_k=\sqrt{\nu}$. Thus from the Lemma \ref{egalitésinus}, we have
$$
\begin{aligned}
\widetilde{T}_2= \displaystyle{\limsup_{k \to +\infty } } \frac{-\log \left \vert\sin\left(\frac{i_k\pi}{\sqrt{\nu}}\right)\right\vert}{i_k^2}&=\displaystyle{\limsup_{k \to +\infty } } \frac{-\log \left \vert\sin\left(\sqrt{\nu}k\pi\right)\right\vert}{\nu k^2}, \\
&=\displaystyle{\limsup_{k \to +\infty } }  \frac{-\log \left \vert\sin\left(\sqrt{\nu}k\pi-h\pi \right)\right\vert}{\nu k^2},\quad \forall h\in \mathbb{N},\\
&\ge\displaystyle{\limsup_{k \to +\infty } }  \frac{-\log \left \vert\sin\left(\sqrt{\nu}q_k\pi-h\pi\right)\right\vert}{\nu q_k^2},\\
&=\displaystyle{\limsup_{k \to +\infty } }  \frac{-\log \left \vert\sin\pi\left(\sqrt{\nu}q_k-p_k\right)\right\vert}{\nu q_k^2},\\
&=\displaystyle{\limsup_{k \to +\infty } }  \frac{-\log \left \vert(\pi\left(\sqrt{\nu}q_k-p_k\right))\right\vert}{\nu q_k^2},\\
&= \displaystyle{\lim_{k \to +\infty } }  \frac{-\log \left \vert\pi q_ke^{-\tau_0 q^2_k}\right\vert}{\nu q_k^2}=\frac{\tau_0}{\nu}.
\end{aligned}
$$
Let us now show $\widetilde{T}_2\le \tau_0$. Note that we can deduce from the previous result that
 
 $$
 \displaystyle{\limsup_{k \to +\infty } }  \frac{-\log \left \vert\sin\pi\left(\sqrt{\nu}q_k-p_k\right)\right\vert}{\nu q_k^2}=\frac{\tau_0}{\nu},
 $$
thus, there exists $ k_0 (\varepsilon) \ge 1 $ such that
 
\begin{equation}
\label{sinuslast}
  \frac{-\log \left \vert\sin\pi\left(\sqrt{\nu}q_k-p_k\right)\right\vert}{\nu q_k^2} \le \frac{\tau_0}{\nu}+ \varepsilon, \quad \forall k\ge k_0(\varepsilon).
\end{equation}
 for every $n\ge 1$, there is $i_n\in \mathbb{N}^*$ for which
 $$
 |\sqrt{\nu}n-i_n|<\frac{1}{2},\quad \forall n\ge 1.
 $$
Let us take $n_0(\varepsilon)=q_{k_0(\varepsilon)} \ge 1$. Thus, using that the sequence $\{ q_k\}_{k \ge 1}$ is strictly increasing, if $n\ge n_0(\varepsilon)$, we infer the existence of $k \ge k_0(\varepsilon)$ such that $q_k \le n < q_{k+1}$. These last properties together with the second formula in (\ref{eq2lem1}) allow us to write
 $$
 |\sqrt{\nu}q_k-p_k| \le |\sqrt{\nu}n-i_n|<\frac{1}{2},
 $$
 and
 $$
 \begin{aligned}
    \vert \sin(\sqrt{\nu}n\pi) \vert   & = \vert \sin(\pi(\sqrt{\nu} n-i_n)) \vert\\
    &=\sin(\pi\vert \sqrt{\nu}n-i_n\vert) , \\
      & \ge  \sin(\pi\vert \sqrt{\nu}q_k-p_k\vert) \\
      &=  \vert \sin(\pi( \sqrt{\nu}q_k-p_k)) \vert,  \quad  \forall n\ge  n_0(\varepsilon), \hbox{ such that }  q_k \le n < q_{k+1}.
 \end{aligned}
 $$
Since  $ q_k \le n < q_{k+1}$, the previous inequality and (\ref{sinuslast}) give
$$
\begin{aligned}
  \frac{-\log \left \vert\sin\left(\sqrt{\nu}n\pi\right)\right\vert}{\nu n^2}&\le   \frac{-\log \left \vert \sin\left(\pi( \sqrt{\nu}q_k-p_k)\right) \right\vert}{\nu q_k^2} \le \frac{\tau_0}{\nu}+\varepsilon, \quad \forall n\ge  n_0(\varepsilon), \hbox{ such that } q_k \le n < q_{k+1},
\end{aligned}
$$
and 
$$
\widetilde{T}_2=\limsup \frac{-\log \left \vert\sin\left(\sqrt{\nu}n\pi\right)\right\vert}{\nu n^2}\le \frac{\tau_0}{\nu}+\varepsilon,\quad \forall \varepsilon>0.
$$
\textbf{Case $\tau_0=\infty$ :} Fix $\sigma>0$ and take $x_0\in (0,+\infty)$ and $\varepsilon>0$ such that  $x_0-\varepsilon >1$. We apply the second item in Lemma \ref{lem1} with,
for instance, $\sigma=1/2$. We deduce the existence of a positive irrational number $\nu >1$ and two sequences of positive integers $\{p_k\}_{k \ge 0}$ and $\{q_k \}_{k\ge 0}$ which fulfills (\ref{eq3lem1}). 
Thus, we deduce

$$
\begin{aligned}
\widetilde{T}_2 &=\displaystyle{\limsup_{k \to +\infty } } \frac{-\log \left \vert\sin\left(\sqrt{\nu}k\pi\right)\right\vert}{\nu k^2}\\
&\ge \displaystyle{\limsup_{k \to +\infty } } \frac{-\log \left \vert\sin\left(\sqrt{\nu} q_k\pi\right)\right\vert}{\nu q_k^2}\\
& = \displaystyle{\limsup_{k \to +\infty } }  \frac{-\log \left \vert\sin\left(\pi\left(q_k\sqrt{\nu}-p_k \right)\right)\right\vert}{\nu q_k^2}\\
&=\displaystyle{\limsup_{k \to +\infty } }  \frac{-\log \left \vert\pi\left(q_k\sqrt{\nu}-p_k \right)\right\vert}{\nu q_k^2}\\
&= \displaystyle{\lim_{k \to +\infty } }  \frac{-\log \left \vert\pi q_ke^{-q_k^{2+\sigma}}\right\vert}{\nu q_k^2}=\infty.
\end{aligned}
$$
This shows that $\widetilde{T}_2 =\infty$ and finishes the third case and the proof of Proposition \ref{prop exmple}. 
$\hfill \square $

\end{document}